\documentclass[12pt]{article}
\usepackage{amsmath,amsthm,amssymb,times,titlesec}
\usepackage{graphicx}
\usepackage{pdfsync}
\numberwithin{equation}{section}

\textheight 
            650pt
\textwidth 
           460pt

\titleformat{\subsubsection}
{\normalfont\fontsize{13}{14}\selectfont\itshape}{\thesubsubsection}{1em}{}

\newcommand{\ee}{\mathrm{e}}

\newcommand{\dt}{\, \mathrm{d}t}
\newcommand{\dx}{\, \mathrm{d}x}
\newcommand{\dy}{\, \mathrm{d}y}

\newcommand{\DD}{{\mathcal D}}

\newcommand{\II}{{\mathcal I}}
\newcommand{\JJ}{{\mathcal J}}

\newcommand{\LL}{{\mathcal L}}

\newcommand{\NN}{{\mathcal N}}

\newcommand{\PP}{{\mathcal P}}

\newcommand{\XX}{{\mathcal X}}

\DeclareMathOperator{\cosech}{cosech}
\DeclareMathOperator{\sech}{sech}
\DeclareMathOperator{\re}{Re}

\newcommand{\sdfrac}[2]{\mbox{\small$\displaystyle\frac{#1}{#2}$}}

\newtheorem{theorem}{Theorem}[section]
\newtheorem{lemma}[theorem]{Lemma}
\newtheorem{proposition}[theorem]{Proposition}

\newtheorem{remark}[theorem]{Remark}
\newtheorem{definition}[theorem]{Definition}
\theoremstyle{definition}

\renewcommand{\i}{\mathrm{i}}

\newcommand\scaleddot{\scalebox{.89}{.}}
\makeatletter
\renewcommand{\dddot}[1]{%
  {\mathop{\kern\z@#1}\limits^{\makebox[0pt][c]{\vbox to-2\ex@{\kern-\tw@\ex@\hbox{\normalfont\scaleddot\kern-0.5pt\scaleddot\kern-0.5pt\scaleddot}\vss}}}}}

\newcounter{count}

\title{Pattern formation on the free surface
of a ferrofluid: spatial dynamics and homoclinic bifurcation}

\author{
M. D. Groves\thanks{Fachrichtung Mathematik, Universit\"{a}t des Saarlandes,
Postfach 151150, 66041 Saarbr\"{u}cken, Germany; 
Department of Mathematical Sciences, Loughborough
University, Loughborough, LE11 3TU, UK
} 
\and D. J. B. Lloyd\thanks{Department of Mathematics, University of Surrey, Guildford, GU2 7XH, UK}
\and A. Stylianou\thanks{Institut f\"{u}r Mathematik, Universit\"{a}t Kassel, 34132 Kassel, Germany}
}
\date{}

\begin{document}

\maketitle

\begin{abstract}
We establish the existence of spatially localised one-dimensional free surfaces of a ferrofluid near onset of
the Rosensweig instability, assuming a general (nonlinear) magnetisation law. It is shown that the ferrohydrostatic equations 
can be derived from a variational principle that allows one to formulate them as an (infinite-dimensional) spatial Hamiltonian system in which the unbounded free-surface direction plays the role of time. A centre-manifold reduction technique converts the problem for small solutions near onset
to an equivalent Hamiltonian system with finitely many degrees of freedom. Normal-form theory yields the
existence of homoclinic solutions to the reduced system, which correspond to spatially localised solutions of the
ferrohydrostatic equations.
\end{abstract}

\section{Introduction}\label{Introduction}

\subsubsection*{Ferrofluids}

The Rosensweig instability - a surface instability of a ferrofluid - has been of interest since the 1960s (see Cowley \& Rosensweig~\cite{CowleyRosensweig67} and Rosensweig \cite{Rosensweig}). In an experiment, a vertical magnetic field is applied to a
static ferrofluid layer, and regular cellular patterns (typically hexagons) emerge on the
fluid surface as the field strength is increased through a critical value
(see Gollwitzer \emph{et al.}~\cite{GollwitzerMatthiesRichterRehbergTobiska07,GollwitzerRehbergRichter10} for recent experimental results). Recently, experiments have shown that spatially localised free-surface structures occur in the hysteresis region between the flat state and the cellular spatially periodic patterns (see Richter~\cite{Richter11},
Richter \& Barashenkov~\cite{RichterBarashenkov05} and Lloyd \emph{et al.}~\cite{LloydGollwitzerRehbergRichter15} for experimental results, and Lavrova \emph{et al.}~\cite{LavrovaMatthiesTobiska08,CaoDing14} for finite-element simulations). Despite this wealth of experimental and numerical evidence little is known theoretically on the existence of localised solutions to the ferrohydrostatic equations.\\

\subsubsection*{Ferrofluids theory}

The first attempt at a theoretical explanation for the nucleation of static cellular patterns was by Gailitis~\cite{Gailitis69,Gailitis77}. Gailitis postulated a free energy for the system (a ferrofluid of infinite depth with a linear magnetisation law)
into which he substituted an \emph{Ansatz} for a cellular free surface. The resulting equations for the unknown coefficients were then solved to find regions of existence for the various cellular patterns (stripes, squares and hexagons). This work was extended by Friedrichs \& Engel~\cite{FriedrichsEngel01}
to finite-depth ferrofluids.

Zaitsev \& Shliomis~\cite{ZaitsevShliomis69} considered one-dimensional static spatially periodic solutions to the 
ferrohydrostatic equations near onset of the Rosensweig instability, again assuming a linear magnetisation law and infinite fluid depth. These results were subsequently extended by Twombly \& Thomas~\cite{TwomblyThomas83} to two-dimensional spatially periodic free-surface solutions and finite-depth fluids. These local bifurcation theories were complemented
by a formal normal-form analysis by Silber \& Knobloch~\cite{SilberKnobloch88} (for two-dimensional patterns and infinite depth),
and there have recently been attempts at deriving time-dependent amplitude equations for spatially periodic free-surface patterns near onset (see Bohlius \emph{et al.}~\cite{BohliusPleinerBrand07,BohilusPleinerBrand11}).

All current theoretical studies, whether rigorous or heuristic,
assume a linear relationship between the magnetisation ${\bf M}$ of the ferrofluid and the
strength of the magnetic field ${\bf H}$. The linearity of the magnetisation law ${\bf M}={\bf M}(|{\bf H}|)$ allows one
to replace a nonlinear equation for the magnetic potential in the fluid by Laplace's equation and greatly simplifies computations
of coefficients in bifurcation and normal-form theory. However, experiments show that realistic magnetisation laws are nonlinear
(e.g.\ see Lloyd \emph{et al.}~\cite{LloydGollwitzerRehbergRichter15} for a discussion of this point), emphasising the need to develop a theoretical framework which includes
general nonlinear magnetisation laws.

\subsubsection*{Localised travelling water waves}

The ferrohydrostatic problem has similarities with the governing equations for travelling gravity-capillary water waves.
One-dimensional spatially localised solutions to that problem (`gravity-capillary solitary waves') with water of finite depth
have been studied intensively in
the last twenty-five years using a technique known as the \emph{Kirchg\"{a}ssner reduction} (see Iooss
\cite{Iooss95},  Groves \& Wahl\'{e}n \cite{GrovesWahlen07} and the references
therein). In this approach one formulates the governing equations as an evolutionary system in which the
horizontal spatial direction plays the role of time (`spatial dynamics'); a centre-manifold reduction principle
is used to show that all spatially localised solutions solve a system of ordinary differential equations, whose solution
set can in principle be determined. A helpful feature of this water-wave
problem is its variational structure, which leads to a Hamiltonian spatial dynamics formulation of the problem;
the Hamiltonian structure is inherited by the reduced system of ordinary differential equations, which can
be treated by well-developed methods for Hamiltonian systems with finitely many degrees of freedom, e.g.\ the
Birkhoff normal-form theory.

\subsubsection*{The present contribution}

\begin{figure}[h]
\centering
\includegraphics[scale=0.8]{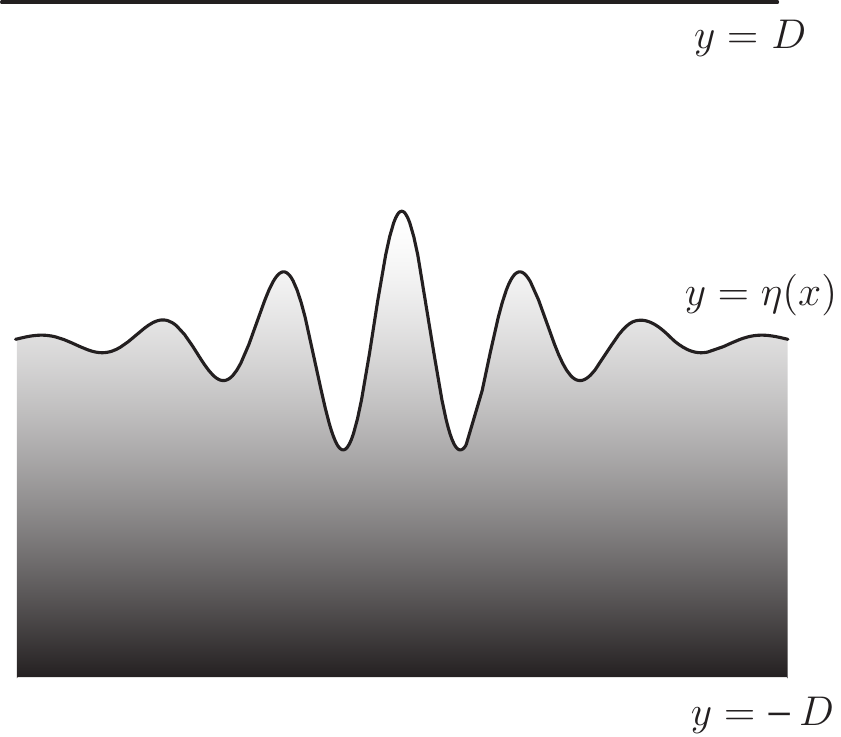}
{\it \caption{Static localised patterns appear on the interface between a non-magnetisable fluid (white) and
a ferrofluid (shaded) as the strength of a vertically aligned magnetic field is varied.}
\label{Definition sketch}}
\end{figure}

We consider two static immiscible perfect fluids in the regions
$$S^\prime:=\{(x,y):\,\eta(x)<y<D\}, \qquad S:=\{(x,y):\,-D<y<\eta(x)\}$$
separated by the free interface $\{y=\eta(x)\}$ (see Figure \ref{Definition sketch}); the upper fluid is non-magnetisable, while the
lower is a ferrofluid with a general nonlinear magnetisation law (see Section \ref{ferroeqns} for a complete mathematical
formulation of the ferrohydrostatic problem). Using spatial dynamics and the Kirchg\"{a}ssner reduction
we present a rigorous existence theory for
small-amplitude localised surface patterns to the ferrohydrostatic problem near onset of the Rosensweig instability.

\begin{figure}[h]
\centering
\includegraphics{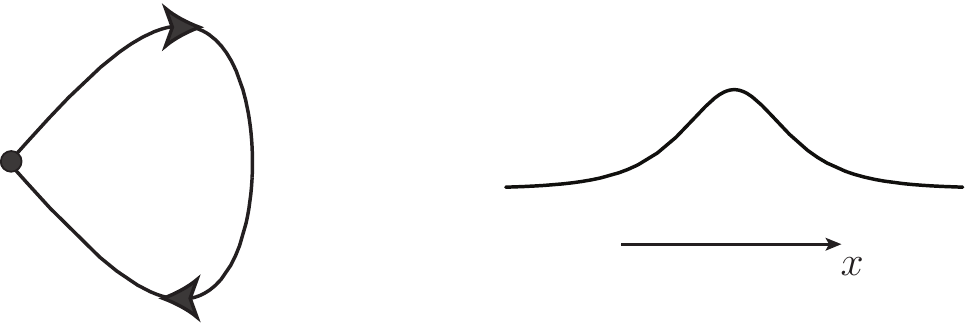}
{\it \caption{A homoclinic solution to the spatial dynamics formulation (left, shown as a trajectory in phase space)
corresponds to a localised solution of the ferrohydrostatic problem (right, sketch of the free surface).}
\label{Homoclinics}}
\end{figure}

Our starting point is the observation that
the governing equations (formulated in terms of
the magnetic potentials $\phi^\prime$ and $\phi$ in respectively the upper and lower fluids)
follow from the variational principle
\begin{align*}
& \delta\left\{-\mu_0\int_{-\infty}^\infty\int_{-D}^{\eta(x)}\big(M(|\nabla \phi|)-M(H)\big)\dy\dx
-\mu_0\int_{-\infty}^\infty\int^{D}_{\eta(x)}\frac{1}{2}|\nabla \phi^\prime|^2\dy\dx\right.  \\
&\qquad\mbox{}+\int_{-\infty}^\infty\Bigg(\frac{1}{2}(\rho-\rho^\prime)g\eta^2
+\mu_0H^2\mu(H)\left(\frac{1}{2}\mu(H)-1\right)\eta\\
&\hspace{1.25in}\mbox{}+\sigma\left(\sqrt{1+\eta_x^2}-1\right) +\mu_0\mu(H)H\left(\phi(-D)-\phi^\prime(D)\right)\Bigg)\Bigg\} \dx=0,
\end{align*}
where $\mu_0$, $g$, $\sigma$ denote respectively the magnetic permeability of free space, acceleration due to gravity and coefficient
of surface tension, $\rho^\prime$ and $\rho$ are the densities of the upper and lower fluids,
$H$ is the strength of the applied magnetic field and
$$\mu(s)=1+\frac{|{\bf M}(s)|}{s}, \qquad M(s) = \int_0^s t \mu(t)\dt;$$
the variations are taken with respect to $\eta$, $\phi^\prime$ and $\phi$ satisfying $\phi^\prime|_{y=\eta}=\phi|_{y=\eta}$.
(Variational (Neumann) boundary conditions at $y=\pm D$ are chosen so that the reference state
$(\eta,\phi^\prime,\phi)=(0,\mu(H)Hy,Hy)$ -- corresponding to a uniform magnetic field and a flat surface --
is a solution to the governing equations.) We regard the above variational functional as an action functional of the form
$$
\int J_\mathrm{f}(\eta, \chi^\prime, \chi, \eta_x, \chi_x^\prime, \chi_x)\dx,
$$
where $\chi^\prime$, $\chi$ are suitably chosen perturbations of the reference states of the magnetic potentials and $J_\mathrm{f}$ is an appropriately transformed version of the `Lagrangian' integrand. Performing a classical Legendre transform yields the desired spatial Hamiltonian formulation
\begin{equation}
\eta_x = \frac{\delta H_\mathrm{f}}{\delta \omega}, \quad
\omega_x = -\frac{\delta H_\mathrm{f}}{\delta \eta}, \quad
\chi^\prime_x = \frac{\delta H_\mathrm{f}}{\delta \xi^\prime}, \quad
\xi^\prime_x = -\frac{\delta H_\mathrm{f}}{\delta \chi^\prime}, \quad
\chi_x = \frac{\delta H_\mathrm{f}}{\delta \xi}, \quad
\xi_x = -\frac{\delta H_\mathrm{f}}{\delta \chi}
\label{Intro - full system}
\end{equation}
of the ferrohydrostatic problem, where $\omega$, $\xi$, $\xi^\prime$ are the momenta associated with the
coordinates $\eta$, $\chi$, $\chi^\prime$ and $H_{\mathrm{f}}$ is the Hamiltonian (see Section \ref{Spatial dynamics}). \emph{Homoclinic solutions} of \eqref{Intro - full system} (solutions with $(\eta,\omega,\chi^\prime,\xi^\prime,\chi,\xi) \rightarrow 0$ as $x \rightarrow \pm \infty$) are of particular interest since they correspond to localised solutions of the ferrohydrostatic problem (see Figure \ref{Homoclinics}).

The Rosensweig instability (the trivial solution becoming unstable to Fourier modes $\pm \ee^{\i s x}$) is associated with a \emph{Hamiltonian Hopf bifurcation} for
\eqref{Intro - full system}: two pairs of simple, purely imaginary
eigenvalues of our linearised Hamiltonian system become complex by colliding on the imaginary axis at the points
$\pm \i s$. Working
in dimensionless variables (see Section \ref{ferroeqns}), one finds that $\pm \i s$ are
double eigenvalues if and only if $s$ is a double root of the equation
$$
s^2\mu(1)(\mu(1)-1)^2  = (\gamma+s^2)\left( \tilde{s}\coth \left(\frac{\tilde{s}}{\beta}\right) + s \mu(1)\coth \left(\frac{s}{\beta}\right)\right), \qquad\hspace{-2.5mm}
\tilde{s}=s\sqrt{\frac{\mu(1)}{\mu(1)+\dot{\mu}(1)}}
$$
(with the assumption $\mu(1)+\dot{\mu}(1)>0$), where
$$\beta=\frac{\sigma}{\mu_0H^2D}, \qquad \gamma=\frac{(\rho-\rho^\prime)\sigma g}{\mu_0^2H^4}.$$
We therefore choose a value $(\beta_0,\gamma_0)$ of $(\beta,\gamma)$ for which this equation has a pair
of double roots and set 
$(\beta,\gamma)=(\beta_0,\gamma_0+\varepsilon)$, so that a Hamiltonian-Hopf bifurcation takes place as $\varepsilon$
increases through zero. In Section \ref{Reduction} we use a centre-manifold reduction principle to
show that our Hamiltonian system near onset admits a locally invariant manifold
of the form
$$\{A e + B f + \overline{A} \overline{e} + \overline{B} \overline{f} + \tilde{r}(A,B,\overline{A},\overline{B},\varepsilon)\},$$
where $e$ and $f$ are suitably normalised generalised eigenvectors at $\varepsilon=0$ and
$\tilde{r}(A,B,\overline{A},\overline{B},\varepsilon)$ satisfies
$\tilde{r} = O(|(A,B,\overline{A},\overline{B})||(\varepsilon,A,B,\overline{A},\overline{B})|)$. The flow on this manifold is described by the
two-degree-of-freedom, reversible Hamiltonian system
\begin{equation}
A_x = \frac{\partial \tilde{H}_\mathrm{f}^\varepsilon}{\partial \overline{B}}, \qquad B_x = -\frac{\partial \tilde{H}_\mathrm{f}^\varepsilon}{\partial \overline{A}}, \label{Intro - reduced system}
\end{equation}
where
$$\tilde{H}_\mathrm{f}^\varepsilon(A,B,\overline{A},\overline{B}) = H_\mathrm{f}^\varepsilon(A e + B f + \overline{A} \overline{e} + \overline{B} \overline{f} + \tilde{r}(A,B,\tilde{A},\overline{B},\varepsilon))$$
(the superscript is added to emphasise the dependence of the Hamiltonian upon $\varepsilon$).
Homoclinic solutions $(A,B)$ to \eqref{Intro - reduced system} (solutions with $A(x)$, $B(x) \rightarrow 0$ as $x \rightarrow \pm\infty$)
generate homoclinic solutions to \eqref{Intro - full system} and hence localised solutions to the ferrohydrostatic
problem.

\begin{figure}[h]
\hspace{0.5cm}\includegraphics[width=4cm]{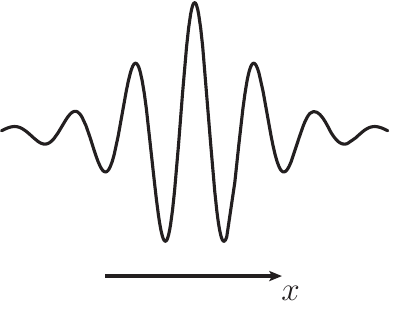}
\hspace{1cm}\includegraphics[width=4cm]{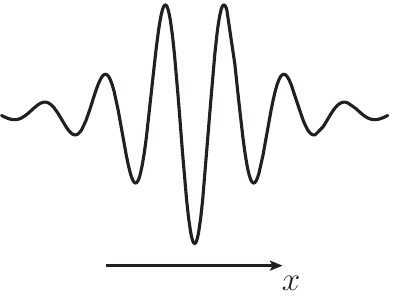}
\hspace{1cm}\includegraphics[width=5cm]{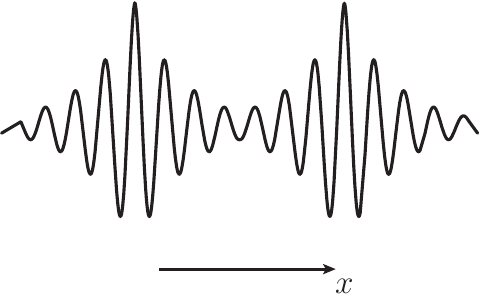}

{\it \caption{Symmetric unipulse patterns (left and centre) co-exist with an infinite
family of multipulse patterns (right).}
\label{Interface sketches}}
\end{figure}

According to the Birkhoff normal-form theory (Section \ref{NF theory}) we can select the coordinates
$A$, $B$, $\overline{A}$, $\overline{B}$ so that the reduced Hamiltonian takes the form
\begin{eqnarray*}
\lefteqn{\tilde{H}_\mathrm{f}^\varepsilon(A,B)=\i \beta_0 q(A\overline{B}-\overline{A}B)+|B|^2}\hspace{0.4in} \nonumber\\
& & \qquad\mbox{} +H_\mathrm{NF}(|A|^2,\i(A\overline{B}-\overline{A}B),\varepsilon)
+ O(|(A,B)|^2|(\varepsilon,A,B)|^{n_0}),
\end{eqnarray*}
where $H_\mathrm{NF}$ is a real polynomial of order $n_0+1$
satisfying
$$H_\mathrm{NF}(|A|^2,\i(A\overline{B}-\overline{A}B),\varepsilon)=O(|(A,B)|^2|(\varepsilon,A,B)|).$$
Existence theories for homoclinic solutions to \eqref{Intro - reduced system} have been given by
Iooss \& P\'{e}rou\`{e}me \cite{IoossPeroueme93} and Buffoni \& Groves \cite{BuffoniGroves99}
under the assumption that the coefficients $c_1$ and $c_3$ in the expansion
\begin{eqnarray*}
\lefteqn{H_\mathrm{NF} = \varepsilon c_1 |A|^2 + \varepsilon \i c_2(A\overline{B}-\overline{A}B)
+c_3|A|^4} \hspace{0.5in} \\
& & \mbox{}+\i c_4|A|^2(A\overline{B}-\overline{A}B)-c_5(A\overline{B}-\overline{A}B)^2
+ \varepsilon^2 c_6 |A|^2 + \varepsilon^2\i c_7 (A\overline{B}-\overline{A}B)+\ldots
\end{eqnarray*}
are respectively negative and positive. Iooss \cite{Iooss95} and Iooss \& P\'{e}rou\`{e}me \cite{IoossPeroueme93} establish the existence of two distinct symmetric homoclinic
solutions with asymptotic expansions
$$A(x) = \pm\left(-\sdfrac{c_3\varepsilon}{c_1}\right)^{\!\!1/2}\!\!\sech \big((-c_1 \varepsilon)^{1/2}x\big)\mathrm{e}^{\mathrm{i}\beta_0 q x}+O(\varepsilon), \qquad B=O(\varepsilon)
$$
as $\varepsilon \rightarrow 0$, while Buffoni \& Groves \cite{BuffoniGroves99} show that
\eqref{Intro - reduced system} 
has an infinite number of geometrically 
distinct homoclinic solutions which generically resemble multiple
copies of one of the Iooss-P\'{e}rou\`{e}me solutions. The free surface is given by the formula
$\eta(x) \sim 2\re A(x) + O(\varepsilon)$, where the constant of proportionality is the first component
of the eigenvector $e$, and Figure \ref{Interface sketches} shows sketches of the free surface corresponding to the
various homoclinic solutions. Explicit formulae for the coefficients $c_1$ and $c_3$ are given in some special cases in
Section \ref{NF theory} (such formulae are unwieldy, and it appears in general more
appropriate to calculate them numerically for a specific choice of $\mu$).

The results presented in this article are mathematically rigorous. In order to make them accessible to as wide a
readership as possible we begin each of Sections \ref{Spatial dynamics} (spatial dynamics),
\ref{Reduction} (centre-manifold reduction) and \ref{NF theory} (homoclinic bifurcation)
with an informal exposition of the theory and a discussion of its computational aspects; the rigorous mathematics
is presented in the second part of each of these sections, specifically Section \ref{FA for SD},
Section \ref{Reduction theory} and Remark \ref{Math of NFT}. The article is intended to be
self-contained, although we omit the details of several lengthy calculations with analogues
in other papers as well as the proofs of specific theorems available elsewhere.

\section{The ferrohydrostatic problem} \label{ferroeqns}
 
We consider two static immiscible perfect fluids in the regions
$$S^\prime:=\{(x,y):\,\eta(x)<y<D\}, \qquad S:=\{(x,y):\,-D<y<\eta(x)\}$$
separated by the free interface $\{y=\eta(x)\}$ (see Figure \ref{Definition sketch}).
The upper, non-magnetisable fluid has unit relative permeability and density $\rho^\prime$,
while the lower is a ferrofluid with density $\rho$.
We denote the magnetic and induction fields in the fluids by respectively
$\mathbf{H^\prime}, \mathbf{H}$ and $\mathbf{B^\prime},\mathbf{B}$, and suppose that the
relationships between them are given by the identities
$$
\mathbf{B^\prime}=\mu_{0}\mathbf{H^\prime}, \qquad
\mathbf{B}	=\mu_{0}(\mathbf{H}+\mathbf{M}(|{\bf H}|)),
$$
where $\mu_{0}$ is the magnetic permeability of free space and $\bf{M}$ is the (prescribed)
magnetic intensity of the ferrofluid. We suppose that ${\bf M}$ and ${\bf H}$ are collinear, so that
$$\mathbf{M}(|{\bf H}|)=|{\bf M}(|{\bf H}|)|\frac{{\bf H}}{|{\bf H}|}.$$
According to Maxwell's equations the magnetic and induction fields are respectively irrotational
and solenoidal, and introducing  magnetic potential functions $\phi^\prime, \phi$ with ${\bf H^\prime}=-\nabla \phi^\prime$,
${\bf H}=-\nabla \phi$, we therefore find that
\begin{equation}
\Delta \phi^\prime=0, \qquad \nabla\cdot(\mu(|\nabla \phi|)\nabla \phi)=0, \label{Basic 1}
\end{equation}
in which
$$\mu(s)=1+\frac{|{\bf M}(s)|}{s}$$
is the magnetic permeability of the ferrofluid relative to that of free space.
(Here, and in the remainder of this paper, equations for `primed' and `non-primed' quantities are supposed to
hold in respectively $S^\prime$ and $S$.) We assume that $\mu$ is a smooth function of $s$ near unity which satisfies
$\mu(1)+\dot{\mu}(1)>0$, where the dot denotes differentiation with respect to $s$.

The ferrohydrostatic Euler equations are given by 
		\begin{align*}
			-\nabla{(p^\prime+\rho^\prime gy)}				& =	\mathbf{0},\\
			\mu_{0}
			\left(\mathbf{M}\cdot
			\nabla\right) {\mathbf{H}}
			-\nabla({p^{\star}+
			\rho gy})					& =	\mathbf{0}
		\end{align*}
(Rosensweig \cite[\S5.1]{Rosensweig}), where $g$ is the acceleration due to gravity, 
$p^\prime$ is the hydrodynamic pressure in the upper fluid and $p^\star$ is the composite pressure in the lower fluid.
The calculation 
$$\left(\mathbf{M}\cdot\nabla\right){\mathbf{H}}=|{\bf M}|\nabla (|{\bf H}|) = \nabla \left( \int_0^{|{\bf H}|} |{\bf M}(t)|\dt\right)$$
shows that these equations are equivalent to
\begin{equation}
-(p^\prime+\rho^\prime gy)=b^\prime_{0}, \qquad \mu_0\int_0^{|{\bf H}|} |{\bf M}(t)|\dt-(p^{\star}+\rho gy) =	b_{0}, \label{Euler}
\end{equation}
where  $b^\prime_{0},b_{0}$ are constants. 

The magnetic boundary conditions at $\{y=\eta(x)\}$ are
$$
\mathbf{H^\prime}\cdot \mathbf{t}=\mathbf{H}\cdot \mathbf{t},
\qquad
\mathbf{B^\prime}\cdot \mathbf{n}=\mathbf{B}\cdot \mathbf{n},		
$$
where 
$$
						\mathbf{t}
			=
						\cfrac{(1,\eta_{{x}})^{\mathrm{T}}
							}{
							\sqrt{1+\eta^{2}_{x}}	},
		\qquad
							\mathbf{n}=\cfrac{(-\eta_x,1)^{\mathrm{T}}
							}{
							\sqrt{1+\eta^{2}_{x}}	}
$$
are the tangent and normal vectors to the interface; 
it follows that 
\begin{equation}
			\phi^\prime-\phi=0,
\qquad	
			\phi^\prime_n
			-\mu(|\nabla\phi|) \phi_n =0
\label{Basic 2}
\end{equation}
for $y=\eta(x)$.
The ferrohydrostatic boundary condition is given by 
		\begin{equation*}
			p^{\star}
			+\frac{\mu_{0}}{2}
			({\mathbf{M}\cdot \mathbf n})^{2}		= p^\prime+2\sigma\kappa
			,	
		\end{equation*} 
(Rosensweig \cite[\S5.2]{Rosensweig}),
in which $\sigma>0$ is the coefficient of surface tension and
$$2\kappa = - \frac{\eta_{xx}}{(1+\eta_x^2)^{3/2}}$$
is the mean curvature of the interface. Using \eqref{Euler}, we find that
$$
\mu_0 \int_0^{|{\bf H}|} |{\bf M}(t)|\dt
+\frac{\mu_{0}}{2}
			({\mathbf{M}\cdot \mathbf n})^{2}
			-(\rho-\rho^\prime)g\eta
			+\frac{\sigma\eta_{xx}}{(1+\eta_x^2)^{3/2}}+b= 0,
$$
where $b=b^\prime_{0}-b_{0}$, or equivalently
\begin{equation}
\mu_0 M(|\nabla \phi|)
-\frac{\mu_0}{2}|\nabla \phi^\prime|^2-\mu_0\big(\mu(|\nabla\phi|)-1\big)\phi_n(\phi_n-\phi_n^\prime)
-(\rho-\rho^\prime)g\eta
+\frac{\sigma\eta_{xx}}{(1+\eta_x^2)^{3/2}}+b	= 0, \label{Basic 3}
\end{equation}
where
$$M(s) = \int_0^s t \mu(t)\dt.$$

The constant $b$ is selected so that 
$(\eta,\phi^\prime,\phi)=(0,\mu(H)Hy,Hy)$ is a solution to \eqref{Basic 1}, \eqref{Basic 2} and \eqref{Basic 3}
(corresponding to a
uniform magnetic field and a flat surface); we therefore set
$b=-\mu_0 M(H)-\mu_0H^2\mu(H)(\frac{1}{2}\mu(H)-1)$. Finally, choosing compatible Neumann boundary conditions
\begin{equation}
\phi_y^\prime|_{y=D}=\mu(H)H,
\qquad
\mu(|\nabla \phi|)\phi_y|_{y=-D}=\mu(H)H.
\label{Basic 4}
\end{equation}
ensures that the system of equations \eqref{Basic 1}, \eqref{Basic 2}--\eqref{Basic 4}
has a variational structure: they follow from the formal variational principle
\begin{align*}
& \delta\left\{-\mu_0\int_{-\infty}^\infty\int_{-D}^{\eta(x)}\big(M(|\nabla \phi|)-M(H)\big)\dy\dx
-\mu_0\int_{-\infty}^\infty\int^{D}_{\eta(x)}\frac{1}{2}|\nabla \phi^\prime|^2\dy\dx\right.  \\
&\qquad\mbox{}+\int_{-\infty}^\infty\Bigg(\frac{1}{2}(\rho-\rho^\prime)g\eta^2
+\mu_0H^2\mu(H)\left(\frac{1}{2}\mu(H)-1\right)\eta\\
&\hspace{1.25in}\mbox{}+\sigma\left(\sqrt{1+\eta_x^2}-1\right) +\mu_0\mu(H)H\left(\phi(-D)-\phi^\prime(D)\right)\Bigg)\Bigg\} \dx=0,
\end{align*}
where the variations are taken with respect to $\eta$, $\phi^\prime$ and $\phi$ satisfying $\phi^\prime|_{y=\eta}=\phi|_{y=\eta}$.

The next step is to introduce dimensionless variables
\[
 (\hat x,\hat y)=\frac{\mu_0H^2}{\sigma}(x,y),\quad \hat \phi:= \frac{\mu_0H}{\sigma}\phi,\quad
 \hat \phi^\prime:= \frac{\mu_0H}{\sigma}\phi^\prime,\quad \hat\eta:= \frac{\mu_0H^2}{\sigma}\eta
\]
and functions
\[
 \hat\mu(s):= \mu(Hs),\qquad \hat M(s):=\frac{1}{H^2}M(Hs)=\int_0^st\hat\mu(t)\dt.
\]
Writing
$(\hat{\eta},\hat{\phi}^\prime,\hat{\phi})=(\hat{\eta},\hat{\psi}^\prime+\mu(1)y,\hat{\psi}+y)$
(so that $(\hat{\eta},\hat{\psi}^\prime,\hat{\psi})=(0,0,0)$ is the `trivial' solution),
we find that
\begin{equation}
\Delta \psi^\prime  = 0, \qquad
\nabla\cdot\left(\mu\left(|\nabla(\psi+y)|\right)\nabla(\psi+y)\right) = 0 \label{Unflattened 1}
\end{equation}
with boundary conditions
\begin{equation}
\psi^\prime_y\big|_{y=\tfrac{1}{\beta}} =0, \qquad
\mu\left(|\nabla(\psi+y)|\right)(\psi_y+1)\big|_{y=-\tfrac{1}{\beta}}-\mu(1) = 0  \label{Unflattened 2}
\end{equation}
and
 \begin{align}
& \psi-\psi^\prime-\left(\rule{0pt}{10pt}\mu(1)-1\right)\eta=0 ,\label{Unflattened 3}\\[1mm]
& \mu\left(|\nabla(\psi+y)|\right)(\psi+y)_n-(\psi^\prime+\mu(1)y)_n=0, \label{Unflattened 4}\\[1mm]
& M\left(|\nabla(\psi+y)|\right)-\frac{1}{2}\left|\nabla\left(\psi^\prime+\mu(1)y\rule{0pt}{10pt}\right)\right|^2
+\sqrt{1+\eta_x^2}\left(\rule{0pt}{10pt} \psi^\prime_y+\mu(1)\right)\left(\rule{0pt}{10pt}\psi^\prime+\mu(1)y\right)_n\nonumber \\
& \quad\mbox{}  -\sqrt{1+\eta_x^2}\ \mu\left(|\nabla(\psi+y)|\rule{0pt}{10pt}\right)\left(\rule{0pt}{10pt}\psi_y+1\right)\left(\rule{0pt}{10pt}\psi+y\right)_n -\gamma\eta \nonumber \\
& \quad\mbox{}  -M(1)-\mu(1)\left(\frac{1}{2}\mu(1)-1\right)+\frac{\eta_{xx}}{(1+\eta_x^2)^{3/2}}=0 \label{Unflattened 5}
 \end{align}
 for $y=\eta(x)$, where
 $$\alpha=\frac{(\rho-\rho^\prime)gD}{\mu_0H^2}, \qquad \beta=\frac{\sigma}{\mu_0H^2D}, \qquad \gamma=\alpha\beta$$
 and the hats have been dropped for notational simplicity.  We use $\beta$ and $\gamma$ as parameters, noting that the
 limit $\beta \rightarrow 0$ corresponds to fluids of infinite depth.
 
Finally, note that equations \eqref{Unflattened 1}--\eqref{Unflattened 5} follow from the formal variational principle $\delta \II=0$, where
\begin{align*}
\II(\eta&,\psi^\prime,\psi):=  \\
& -\int_{-\infty}^\infty\int_{-\tfrac{1}{\beta}}^{\eta(x)}\big(M(|\nabla\left( \psi+y\right)|)-M(1)\big)\dy\dx
-\int_{-\infty}^\infty\int^{\tfrac{1}{\beta}}_{\eta(x)}\frac{1}{2}|\nabla\left( \psi^\prime+\mu(1)y\right)|^2\dy\dx  \\
&\quad\mbox{}+\int_{-\infty}^\infty\Bigg\{\frac{\gamma}{2}\eta^2
+\mu(1)\left(\frac{1}{2}\mu(1)-1\right)\eta
+\left(\sqrt{1+\eta_x^2}-1\right)
+\mu(1)\left(\psi(-\tfrac{1}{\beta})-\psi^\prime(\tfrac{1}{\beta})\right)\Bigg\} \dx 
\end{align*}
and the variations are taken with respect to $\eta$, $\psi^\prime$ and $\psi$ with
$\psi^\prime(0)=\psi(0)-(\mu(1)-1)\eta$.

\section{Spatial dynamics} \label{Spatial dynamics}

\subsection{Formulation as a spatial Hamiltonian system} \label{SD - derivation}
The first step is to use the `flattening' transformation
$$
Y=\left\{
\begin{aligned}
 \frac{y-\eta(x)}{1-\beta\eta(x)}, & \qquad y \geq \eta(x), \\
 \frac{y-\eta(x)}{1+\beta\eta(x)}, & \qquad y < \eta(x),
\end{aligned}
\right.
$$
to map the variable fluid domains $S^\prime$ and $S$ into fixed strips ${\mathbb R} \times (0,\tfrac{1}{\beta})$ and ${\mathbb R} \times (-\tfrac{1}{\beta},0)$ and the free interface $\{y=\eta(x)\}$ into $\{y=0\}$.
Replacing the symbol $Y$ by $y$ for notational simplicity, we find that the corresponding `flattened' variables
$$\chi^\prime(x,Y)=\psi^\prime(x,y), \qquad \chi(x,Y)=\psi(x,y)$$
satisfy the equations
\begin{align}
   \chi^\prime_{xx}-K^\prime_1\eta_{xx}\chi^\prime_y-2K^\prime_1\eta_x\chi^\prime_{xy}+K_1^{\prime 2}\eta_x^2 \chi^\prime_{yy}
  +K_2^{\prime 2}\chi^\prime_{yy}-2\beta K^\prime_1 K^\prime_2\eta_x^2\chi^\prime_y= & 0, \label{Flattened 1} \\[2mm]
    \mu^\star\left(
\chi_{xx}-K_1\eta_{xx}\chi_y-2K_1\eta_x\chi_{xy}+K_1^2\eta_x^2 \chi_{yy}
  +K_2^2\chi_{yy}+2\beta K_1 K_2\eta_x^2\chi_y \right)\hspace{0.5cm}{}& \nonumber\\
  +\left(\chi_x-\eta_xK_1\chi_y,K_2\chi_y+1\right)\cdot (\mu^\star_x-\eta_x K_1\mu^\star_y,K_2\mu^\star_y)= & 0 \label{Flattened 2}
  \end{align}
  with boundary conditions
\begin{equation}
   \mu^\star(K_2\chi_y+1)\big|_{y=-\tfrac{1}{\beta}}-\mu(1) = 0, \qquad
      K^\prime_2\chi_y^\prime\big|_{y=\tfrac{1}{\beta}} = 0 \label{Flattened 3}
   \end{equation}
and
{\allowdisplaybreaks\begin{align}
&  \chi-\chi^\prime-\left(\rule{0pt}{10pt}\mu(1)-1\right)\eta=0, \label{Flattened 4}\\
&  \left(1+\eta_x^2\right)(\mu^\star K_2\chi_y-K^\prime_2\chi^\prime_{y})
  -\eta_x(\mu^\star\chi_x-\chi^\prime_x)+\mu^\star-\mu(1)=0, \label{Flattened 5}\\[2mm]
&  M^\star
  -\frac{1}{2}\chi_x^{\prime 2}+\left(1+\eta_x^2\right)\left(\mu(1)K^\prime_2\chi^\prime_y-\mu^\star K_2\chi_y\right) \nonumber \\
& \quad\mbox{}  +\left(1+\eta_x^2\right)\left(\frac{1}{2}\left(K^\prime_2\chi^\prime_y\right)^2-\mu^\star\left(K_2\chi_y\right)^2\right)
  -\eta_x(\mu(1)\chi_x^\prime-\mu^\star\chi_x) \nonumber \\
& \quad\mbox{}   +\mu^\star\left(K_2\chi_y\eta_x\chi_x-K_2\chi_y-1\right)
 -\gamma\eta-M(1)+\mu(1)+\frac{\eta_{xx}}{(1+\eta_x^2)^{3/2}}=0 \label{Flattened 6}
 \end{align}}
on $y=0$; here
\begin{align*}
\mu^\star & = \mu \left(\sqrt{(\chi_x-\eta_xK_1\chi_y)^2
 +(K_2\chi_y+1)^2}\right), \\
M^\star &= M \left(\sqrt{(\chi_x-\eta_xK_1\chi_y)^2
 +(K_2\chi_y+1)^2}\right)
\end{align*}
and
$$K^\prime_1= \frac{1-\beta y}{1-\beta \eta}, \quad K_1= \frac{1+\beta y}{1+\beta \eta}, \qquad
K^\prime_2 = \frac{1}{1-\beta \eta}, \quad K_2=\frac{1}{1+\beta\eta}.$$

Observe that equations \eqref{Flattened 1}--\eqref{Flattened 6} follow from the new variational principle
$\delta \JJ=0$, where
\begin{align}
  \JJ(&\eta,\chi^\prime,\chi):=
  -\int_{-\infty}^\infty\int_{-\tfrac{1}{\beta}}^0 \frac{1}{K_2} \big(M^\star-M(1)\big) \dy\dx \nonumber\\
&\mbox{}  -\int_{-\infty}^\infty\int_{0}^{\tfrac{1}{\beta}} \frac{1}{2K^\prime_2} \left\{\rule{-3pt}{12pt}\left(\chi^\prime_x-K^\prime_1\eta_x\chi^\prime_y\right)^2+\left(K^\prime_2\chi^\prime_y+\mu(1)\right)^2\right\} \dy\dx \nonumber \\
&\mbox{}+\int_{-\infty}^\infty\Bigg\{\frac{\gamma}{2}\eta^2
+\mu(1)\left(\frac{1}{2}\mu(1)-1\right)\eta
+\left(\sqrt{1+\eta_x^2}-1\right)
+\mu(1)\left(\psi(-\tfrac{1}{\beta})-\psi^\prime(\tfrac{1}{\beta})\right)\Bigg\} \dx \label{Definition of cal J}
 \end{align}
 and the variations are taken in $\eta$, $\chi^\prime$ and $\chi$ satisfying the side constraint
 \begin{equation}
 \chi^\prime(0)=\chi(0)-(\mu(1)-1)\eta \label{Side constraint}
 \end{equation}
 (the functional $\JJ$ is obtained from $\II$ by
 `flattening').
We exploit this variational principle by regarding $\JJ$ as an action functional
of the form
$$\JJ = \int J_\mathrm{f}(\eta, \chi^\prime, \chi, \eta_x, \chi_x^\prime, \chi_x) \dx,$$
in which $J$ is the integrand on the right-hand side of equation \eqref{Definition of cal J},
and deriving a canonical Hamiltonian formulation of \eqref{Flattened 1}--\eqref{Flattened 6} by means of the Legendre transform.
To this end, let us introduce new variables $\omega$, $\xi^\prime$ and $\xi$ by the formulae
\begin{align*}
\omega &= \frac{\delta J_\mathrm{f}}{\delta \eta_x}
=\int_{-\tfrac{1}{\beta}}^0 \frac{K_1}{K_2}\mu^\star(\chi_x-K_1\eta_x\chi_y)\chi_y\dy
+\int_0^{\tfrac{1}{\beta}} \frac{K^\prime_1}{K^\prime_2}(\chi_x^\prime-K^\prime_1\eta_x\chi_y^\prime)\chi_y^\prime\dy
+\frac{\eta_x}{\sqrt{1+\eta_x^2}}, \\
\xi^\prime &= \frac{\delta J_\mathrm{f}}{\delta \chi_x^\prime} = -\frac{1}{K^\prime_2}(\chi_x^\prime-K^\prime_1\eta_x\chi_y^\prime), \\
\xi &= \frac{\delta J_\mathrm{f}}{\delta \chi_x} = -\frac{1}{K_2}\mu^\star(\chi_x-K_1\eta_x\chi_y)
\end{align*}
and define the Hamiltonian function by
\begin{align}
H_\mathrm{f}(\eta,\omega,\chi^\prime,\xi^\prime,\chi,\xi)
& = \int_{-\tfrac{1}{\beta}}^0 \xi \chi_y \dy + \int_0^{\tfrac{1}{\beta}} \xi^\prime \chi_x^\prime \dy +\omega\eta_x-J(\eta, \chi^\prime, \chi, \eta_x, \chi_x^\prime, \chi_x) \nonumber \\
& = \int_0^{\tfrac{1}{\beta}} \frac{1}{2}K^\prime_2(\chi_y^{\prime 2}-\xi^{\prime 2})\dy
+\int_{-\tfrac{1}{\beta}}^0 \left\{ \frac{1}{K_2}\big(M^\dagger-M(1)\big)-\frac{K_2}{\mu^\dagger}\xi^2\right\}\dy \nonumber \\
& \qquad\mbox{}-\mu(1)(\chi(0)+\chi(-\tfrac{1}{\beta}))-\frac{\gamma}{2}\eta^2+1-\sqrt{1-W^2}, \label{Hamiltonian}
\end{align}
in which
\begin{align*}
W &=\omega + \int_{-\tfrac{1}{\beta}}^0 K_1 \xi\chi_y\dy + \int_0^{\tfrac{1}{\beta}} K^\prime_1 \xi^\prime \chi_y^\prime \dy, \\
\mu^\dagger & = \mu\left(\sqrt{\nu^{-1}(K_2\xi;K_2\chi_y+1)^2+(K_2\chi_y+1)^2}\right), \\
M^\dagger & = M\left(\sqrt{\nu^{-1}(K_2\xi;K_2\chi_y+1)^2+(K_2\chi_y+1)^2}\right)
\end{align*}
and $\nu(\cdot;t): {\mathbb R} \rightarrow {\mathbb R}$ is given by the formula $\nu(s;t)=s\mu(\sqrt{s^2+t^2})$, where $t$ is
a parameter whose value is near unity. (Using the the calculations $\nu(0,1)=0$ and $\partial_1\nu(0,1)=\mu(1)>1$, one
finds from the inverse-function theorem that $\nu$ is invertible for $(s,t)$ near $(0,1)$ (see Section \ref{FA for SD}).)

Hamilton's equations are given explicitly by
{\allowdisplaybreaks\begin{align}
\eta_x = \frac{\delta H_\mathrm{f}}{\delta \omega} &= \frac{W}{\sqrt{1-W^2}}, \label{First HE 1}\\
\omega_x = -\frac{\delta H_\mathrm{f}}{\delta \eta} & = 
-\beta\int_{-\tfrac{1}{\beta}}^0 \left( M^\dagger-M(1)-\mu^\dagger K_2\chi_y (K_2\chi_y+1)\right)\dy
+ \beta\int_0^{\tfrac{1}{\beta}} \frac{1}{2}K_2^{\prime 2}(\xi^{\prime 2}-\chi_y^{\prime 2})\dy
\nonumber \\
& \qquad\mbox{}
+\frac{\beta W}{\sqrt{1-W^2}}\left(\int_{-\tfrac{1}{\beta}}^0 K_1K_2\xi\chi_y \dy
-\int_0^{\tfrac{1}{\beta}} K^\prime_1K^\prime_2\xi^\prime\chi^\prime_y \dy\right) \nonumber \\
& \qquad \mbox{}-\frac{W}{\sqrt{1-W^2}}(\mu(1)K_2^\prime\xi^\prime(0)-K_2\xi(0)) \nonumber \\
& \qquad\mbox{}+\mu^\dagger|_{y=0}(K_2\chi_y(0)+1)-\mu(1)(K_2^\prime\chi_y^\prime(0)+1)
+\gamma\eta, \label{First HE 2}\\
\chi^\prime_x = \frac{\delta H_\mathrm{f}}{\delta \xi^\prime} &= -K^\prime_2\xi^\prime + \frac{WK^\prime_1\chi_y^\prime}{\sqrt{1-W^2}}, \label{First HE 3}\\
\xi^\prime_x = -\frac{\delta H_\mathrm{f}}{\delta \chi^\prime} &=K^\prime_2\chi_{yy}^\prime + \frac{W(K_1^\prime\xi^\prime)_y}{\sqrt{1-W^2}}, \label{First HE 4}\\
\chi_x = \frac{\delta H_\mathrm{f}}{\delta \xi} &=-\frac{K_2}{\mu^\dagger}\xi+ \frac{WK_1\chi_y}{\sqrt{1-W^2}}, \label{First HE 5}\\
\xi_x = -\frac{\delta H_\mathrm{f}}{\delta \chi} &=((\mu^\dagger(K_2\chi_y+1))_y + \frac{W(K_1\xi)_y}{\sqrt{1-W^2}} \label{First HE 6}
\end{align}}
and are accompanied by the side constraint \eqref{Side constraint} and further boundary conditions
\begin{align}
& \chi^\prime_y(\tfrac{1}{\beta})=0, \label{BC 1}\\
& \mu^\dagger|_{y=-\tfrac{1}{\beta}} \left(K_2\chi_y(-\tfrac{1}{\beta})+1\right)-\mu(1)=0, \label{BC 2}\\
& \displaystyle \frac{K_2W\xi(0)}{\sqrt{1-W^2}}+\mu^\dagger|_{y=0} \left(K_2\chi_y(0)+1\right)-\mu(1)-\frac{K^\prime_2W \xi^\prime(0) }{\sqrt{1-W^2}}-K^\prime_2\chi^\prime_y(0)=0, \label{BC 3}\\
&\displaystyle \frac{1}{\mu^\dagger|_{y=0}}K_2 \xi(0)-\frac{K_2 W \chi_y(0)}{\sqrt{1-W^2}}-K^\prime_2 \xi^\prime(0)+\frac{K^\prime_2 W \chi^\prime_y(0)}{\sqrt{1-W^2}}+\left(\mu(1)-1\right)\frac{W}{\sqrt{1-W^2}}=0. \label{BC 4}
\end{align}
The first three of these equations arise from the integration
by parts necessary to compute \eqref{First HE 4} and \eqref{First HE 6} subject to the constraint
\eqref{Side constraint}, while the fourth is the compatibility condition
which ensures that $\chi_x^\prime(0)=\chi_x(0)-(\mu(1)-1)\eta_x$.
Note further that our equations are reversible, that is invariant under the transformation
$(\eta,\omega,\chi^\prime,\xi^\prime,\chi,\xi)(x) \mapsto R(\eta,\omega,\chi^\prime,\xi^\prime,\chi,\xi)(-x)$, where
the \emph{reverser} is
defined by $R(\eta,\omega,\chi^\prime,\xi^\prime,\chi,\xi)=(\eta,-\omega,\chi^\prime,-\xi^\prime,\chi,-\xi)$.

Our equations are also invariant under the transformation
$\chi^\prime \mapsto \chi^\prime + c$, $\chi \mapsto \chi + c$ for any constant $c$. To eliminate this symmetry it is
convenient to replace $(\xi^\prime,\chi^\prime,\xi,\chi)$ with new variables $(\overline{\chi}^\prime,\overline{\xi}^\prime,
\overline{\chi},\overline{\xi},\hat{\chi},\hat{\xi})$, where
$$\overline{\chi}^\prime:= \chi^\prime - \hat{\chi}, \quad \overline{\chi}= \chi - \hat{\chi}, \qquad
\overline{\xi}^\prime:= \xi^\prime - \hat{\xi}, \quad \overline{\xi}= \xi - \hat{\xi},$$
and
$$
\hat{\chi} := \frac{1}{2}\left(\int_0^{\tfrac{1}{\beta}} \chi^\prime \dy +\int_{-\tfrac{1}{\beta}}^0 \chi \dy\right),
\qquad
\hat{\xi} := \frac{1}{2}\left(\int_0^{\tfrac{1}{\beta}} \xi^\prime \dy +\int_{-\tfrac{1}{\beta}}^0 \xi \dy\right).
$$
This transformation leads to a new canonical Hamiltonian system with Hamiltonian
\begin{align*}
\overline{H}_\mathrm{f}(\eta,\omega,\overline{\chi}^\prime,\overline{\xi}^\prime,\hat{\chi}^\prime,\hat{\xi}^\prime,\overline{\chi},\overline{\xi},\hat{\chi},\hat{\xi})
& =H_\mathrm{f}(\eta,\omega,\overline{\chi}^\prime+\hat{\chi},\overline{\xi}^\prime+\hat{\xi},
\overline{\chi}+\hat{\chi},\overline{\xi}+\hat{\xi}) \\
& = H_\mathrm{f}(\eta,\omega,\overline{\chi}^\prime,\overline{\xi}^\prime+\hat{\xi},
\overline{\chi},\overline{\xi}+\hat{\xi})
\end{align*}
and additional constraints
\begin{equation}
\int_0^{\tfrac{1}{\beta}} \overline{\chi}^\prime \dy + \int_{-\tfrac{1}{\beta}}^0 \overline{\chi} \dy =0, 
\qquad \int_0^{\tfrac{1}{\beta}} \overline{\xi}^\prime \dy + \int_{-\tfrac{1}{\beta}}^0 \overline{\xi} \dy =0.
\label{Integral constraints}
\end{equation}
Observe that $\hat{\chi}$ is a cyclic variable whose conjugate $\hat{\xi}$ is a conserved quantity; we proceed in standard fashion
by setting $\hat{\xi}=0$,
considering the equations for $(\eta,\omega,\overline{\chi}^\prime,\overline{\xi}^\prime,\overline{\chi},\overline{\xi})$,
and recovering $\hat{\chi}$ by quadrature. Dropping the bars for notational simplicity, one finds that
Hamilton's equations for the reduced system are
{\allowdisplaybreaks\begin{align}
\eta_x &= \frac{W}{\sqrt{1-W^2}}, \label{HE 1} \\
\omega_x & = 
-\beta\int_{-\tfrac{1}{\beta}}^0 \left( M^\dagger-M(1)-\mu^\dagger K_2\chi_y (K_2\chi_y+1)\right)\dy
+\beta\int_0^{\tfrac{1}{\beta}} \frac{1}{2}K_2^{\prime 2}(\xi^{\prime 2}-\chi_y^{\prime 2})\dy
\nonumber \\
& \qquad\mbox{}
-\frac{\beta W}{\sqrt{1-W^2}}\left(\int_{-\tfrac{1}{\beta}}^0 K_1K_2\xi\chi_y \dy
+\int_0^{\tfrac{1}{\beta}} K^\prime_1K^\prime_2\xi^\prime\chi^\prime_y \dy\right) \nonumber \\
& \qquad \mbox{}+\frac{W}{\sqrt{1-W^2}}(\mu(1)K_2^\prime\xi^\prime(0)-K_2\xi(0)) \nonumber \\
& \qquad\mbox{}+\mu^\dagger|_{y=0}(K_2\chi_y(0)+1)-\mu(1)(K_2^\prime\chi_y^\prime(0)+1)
+\gamma\eta, \label{HE 2}\\
\chi_x^\prime &= -K^\prime_2\xi^\prime + \frac{WK^\prime_1\chi_y^\prime}{\sqrt{1-W^2}} \nonumber\\
& \qquad \mbox{}
-\frac{1}{2}\left(\int_0^{\tfrac{1}{\beta}}\!\! \left(-K^\prime_2\xi^\prime + \frac{WK^\prime_1\chi_y^\prime}{\sqrt{1-W^2}}\right)\!\!\dy
+\int_{-\tfrac{1}{\beta}}^0 \!\!\left(-\frac{K_2}{\mu^\dagger}\xi
+ \frac{WK_1\chi_y}{\sqrt{1-W^2}}\right)\!\!\dy\right)\!, \label{HE 3}\\
\xi_x^\prime &=K^\prime_2\chi_{yy}^\prime + \frac{W(K_1\xi)_y}{\sqrt{1-W^2}}, \label{HE 4} \\
\chi_x &=-\frac{K_2}{\mu^\dagger}\xi
+ \frac{WK_1\chi_y}{\sqrt{1-W^2}} \nonumber \\
& \qquad \mbox{}
-\frac{1}{2}\left(\int_0^{\tfrac{1}{\beta}}\!\! \left(-K^\prime_2\xi^\prime + \frac{WK^\prime_1\chi_y^\prime}{\sqrt{1-W^2}}\right)\!\!\dy
+\int_{-\tfrac{1}{\beta}}^0\!\! \left(-\frac{K_2}{\mu^\dagger}\xi
+ \frac{WK_1\chi_y}{\sqrt{1-W^2}}\right)\!\!\dy\right)\!, \label{HE 5} \\
\xi_x &=((\mu^\dagger(K_2\chi_y+1))_y + \frac{W(K_1\xi)_y}{\sqrt{1-W^2}}, \label{HE 6}
\end{align}}
with constraints \eqref{Side constraint}, \eqref{Integral constraints} and boundary conditions \eqref{BC 1}--\eqref{BC 4};
the quantity $\hat{\chi}$ is recovered by quadrature from the equation
\begin{equation}
\hat{\chi}_x =
\frac{1}{2}\left(\int_0^{\tfrac{1}{\beta}}\!\! \left(-K^\prime_2\xi^\prime + \frac{WK^\prime_1\chi_y^\prime}{\sqrt{1-W^2}}\right)\!\!\dy
+\int_{-\tfrac{1}{\beta}}^0 \!\!\left(-\frac{K_2}{\mu^\dagger}\xi
+ \frac{WK_1\chi_y}{\sqrt{1-W^2}}\right)\!\!\dy\right). \label{Quadrature}
\end{equation}
Notice that these equations are also reversible (with respect to the same reverser $R)$.

\subsection{Functional-analytic basis} \label{FA for SD} 

We begin with a precise statement of the invertibility of the function
$\nu(\cdot;t): {\mathbb R} \rightarrow {\mathbb R}$ defined by
the formula $\nu(s;t)=s\mu(\sqrt{s^2+t^2})$.

\begin{proposition} \label{Invert nu}
There exist open neighbourhoods $\Lambda_1$, $\Lambda_3$ of the origin and $\Lambda_2$ of unity in ${\mathbb R}$ such that
$\nu(\cdot,t): \Lambda_1 \rightarrow \Lambda_3$ is a bijection for each $t \in \Lambda_2$. Furthermore
$\nu \in C^\infty(\Lambda_1 \times \Lambda_2)$ and $\nu^{-1} \in C^\infty(\Lambda_1 \times \Lambda_3)$.
\end{proposition}

Next we recall the differential-geometric definitions of a Hamiltonian system and Hamilton's equations for its associated
vector field.

\begin{definition}
A \underline{Hamiltonian system} consists of
a triple $(M,\Omega,H)$, where $M$ is a manifold, $\Omega: TM \times TM \rightarrow {\mathbb R}$
is a closed, weakly nondegenerate bilinear form (the \underline{symplectic $2$-form}) and the
\underline{Hamiltonian}
$H: N \rightarrow {\mathbb R}$ is a smooth function on a manifold domain $N$ of $M$ (that is, a
manifold $N$ which is smoothly embedded in $M$ and has the property that $TN|_{n}$ is densely embedded
in $TM|_{n}$ for each $n \in N$).
\end{definition}

To apply this definition to the Hamiltonian system derived in Section \ref{SD - derivation} above, we introduce the Hilbert spaces
\begin{align*}
M_s & = \Bigg\{(\eta,\omega,\chi^\prime,\xi^\prime,\chi,\xi)\in {\mathbb R} \times {\mathbb R} \times H^{s+1}(0,\tfrac{1}{\beta}) \times H^s(0,\tfrac{1}{\beta}) \times H^{s+1}(-\tfrac{1}{\beta},0) \times H^s(-\tfrac{1}{\beta},0)\!: \\[-2mm]
& \hspace{1in} \chi^\prime(0)=\chi(0)-(\mu(1)-1)\eta, \\
& \hspace{1in}\int_0^{\tfrac{1}{\beta}} \chi^\prime \dy + \int_{-\tfrac{1}{\beta}}^0 \chi \dy =0, 
\ \int_0^{\tfrac{1}{\beta}} \xi^\prime \dy + \int_{-\tfrac{1}{\beta}}^0 \xi \dy =0\Bigg\},
 \qquad s \geq 0,
\end{align*}
and let $N_1$ be a neighbourhood of the origin in $M_1$ such that $|W|<1$, $|\eta|<\tfrac{1}{\beta}$ and
$K_2\chi_y(y)+1 \in \Lambda_2$, $K_2\xi(y) \in \Lambda_3$ for each $y \in [0,1]$ (recall that $H^1(-\frac{1}{\beta},0)$ is
continuously embedded in $C[-\frac{1}{\beta},0]$).
Observe that $N_1$ is a manifold domain of $M_0$, while the formula
\begin{align*}
\Omega_\mathrm{f}((\eta_1,\omega_1,&\chi^\prime_1,\xi^\prime_1,\chi_1,\xi_1), (\eta_2,\omega_2,\chi^\prime_2,\xi^\prime_2,\chi_2,\xi_2)) \\
& =\omega_2\eta_1-\eta_2\omega_1 + \int_0^{\tfrac{1}{\beta}} (\xi^\prime_2\chi^\prime_1-\chi_2^\prime\xi^\prime_1) \dy
 + \int_{-\tfrac{1}{\beta}}^0 (\xi_2\chi_1-\chi_2\xi_1) \dy,
 \end{align*}
defines a weakly nondegenerate bilinear form $M_0 \times M_0 \rightarrow {\mathbb R}$ and hence a constant
symplectic $2$-form $TM_0 \times TM_0 \rightarrow {\mathbb R}$ (its closure follows from the fact that it is constant).
Furthermore, the function $H_\mathrm{f}$ given by \eqref{Hamiltonian} belongs to
$C^\infty(N_1,{\mathbb R})$, so that the triple $(M_0,\Omega_\mathrm{f},H_\mathrm{f})$ is a Hamiltonian system.

\begin{definition}
Consider a Hamiltonian system $(M,\Omega,H)$, where $H \in C^\infty(N,{\mathbb R})$ and $N$ is a manifold
domain of $M$. Its \underline{Hamiltonian vector field} $v_H$ with domain $\DD(v_H) \subseteq N$ is defined
as follows. The point $n \in N$ belongs to $\DD(v_H)$ with $v_H|_{n} := \hat{v}\in TM|_{n}$ if and only if
$$\Omega|_{n}(\hat{v}, v) = {\bf d}H|_{n}(v)$$
for all tangent vectors $v \in TM|_{n}$ (by construction ${\bf d}H|_{n} \in T^\ast N|_{n}$ admits a unique extension
${\bf d}H|_{n} \in T^\star M|_{n}$.) \underline{Hamilton's equations} for $(M,\Omega,H)$ are the differential equations
$$u_x = v_H(u):= v_H|_u$$
which determine the trajectories $u \in C^1({\mathbb R},M_0) \cap C({\mathbb R},N_1)$ of its Hamiltonian vector field.
\end{definition}

Applying the criterion in this definition to the ferrofluid Hamiltonian system $(M_0,\Omega_\mathrm{f},H_\mathrm{f})$, one finds that
$$\DD(v_{H_\mathrm{f}}) =\{(\eta,\omega,\chi^\prime,\xi^\prime,\chi,\xi) \in N_1: B(\eta,\omega,\chi^\prime,\xi^\prime,\chi,\xi)=0\},$$
where $B \in C^\infty(N_1, {\mathbb R}^4)$ is defined by the left-hand sides of equations \eqref{BC 1}--\eqref{BC 4},
and that Hamilton's equations
are given explicitly by \eqref{First HE 1}--\eqref{First HE 6}. Observe that these equations constitute a
quasilinear evolutionary system in the phase space $M_0$ with nonlinear boundary conditions
\eqref{BC 1}--\eqref{BC 4}; its right-hand side is a smooth mapping $N_1 \rightarrow M_0$.

It remains to confirm the relationship between a solution to Hamilton's equations for\linebreak
$(M_0,\Omega_\mathrm{f},H_\mathrm{f})$
and a solution to the `flattened' ferrohydrostatic problem \eqref{Flattened 1}--\eqref{Flattened 6}. 
Suppose that $(\eta,\omega,\chi^\prime,\xi^\prime,\chi,\xi)$ is a smooth
solution of Hamilton's equations and compute $\hat{\chi}$ from
\eqref{Quadrature} by quadrature. An explicit calculation shows that the variables
$\tilde{\eta}$, $\tilde{\chi}^\prime$, $\tilde{\chi}$ given by $\tilde{\eta}(x)=\eta(x)$,
$\tilde{\chi}^\prime(x,z) = \chi^\prime(x)(z)+\hat{\chi}(x)(z)$, $\tilde{\chi}(x,z) = \chi(x)(z)+\hat{\chi}(x)(z)$ solve 
\eqref{Flattened 1}--\eqref{Flattened 6} (see Buffoni, Groves \& Toland \cite[Theorem 2.1]{BuffoniGrovesToland96} for a
discussion of this procedure in the context of water waves).

\section{Centre-manifold reduction} \label{Reduction}

\subsection{Reduction to a two-degree of freedom Hamiltonian system}

We now set $(\beta,\gamma)=(\beta_0,\gamma_0+\varepsilon)$, where the values of $\beta_0$, $\gamma_0$
are appropriately chosen and fixed (see below) and $\varepsilon$ plays the role of a bifurcation parameter.
The Hamiltonian formulation of our ferrohydrostatic problem is accordingly written as
\begin{equation}
u_x = f^\varepsilon(u), \label{HE}
\end{equation}
where $u=(\eta,\omega,\chi^\prime,\xi^\prime,\chi,\xi)$ and $f^\varepsilon$ is given by the right-hand side of \eqref{HE 1}--\eqref{HE 6}, with linear constraints
\eqref{Side constraint}, \eqref{Integral constraints} and nonlinear boundary conditions
\begin{equation}
B(u)=0, \label{NLBC}
\end{equation}
where $B$ is given by the left-hand sides of equations \eqref{BC 1}--\eqref{BC 4}. Similarly, we denote the Hamiltonian
\eqref{Hamiltonian} with this parameter choice by $H^\varepsilon_\mathrm{f}$.

The corresponding linearised system is
$$u_x = Lu,$$
where

\begin{equation}
L\begin{pmatrix} \eta \\ \omega \\ \tau^\prime \\ \zeta^\prime \\ \tau \\ \zeta \end{pmatrix}
=
\begin{pmatrix}
\omega\\[0.5em]
\displaystyle \left(\mu(1)+\dot{\mu}(1)\right)\tau_y(0)-\mu(1)\tau^\prime_y(0) +\gamma_0\eta\\[0.5em]
\displaystyle -\zeta^\prime
+\frac{1}{2\mu(1)}\int_{-\tfrac{1}{\beta_0}}^0 \zeta \dy + \frac{1}{2}\int_{0}^{\tfrac{1}{\beta_0}}\zeta^\prime\dy\\[1.5em]
\displaystyle  \tau^\prime_{yy}\\[0.5em]
\displaystyle -\frac{1}{\mu(1)}\zeta+\frac{1}{2\mu(1)}\int_{-\tfrac{1}{\beta_0}}^0 \zeta \dy + \frac{1}{2}\int_{0}^{\tfrac{1}{\beta_0}}\zeta^\prime\dy\\[1.5em]
\displaystyle (\mu(1)+\dot{\mu}(1))\tau_{yy}
\end{pmatrix},
\label{Explicit formula for L}
\end{equation}\\ \\
with linear constraints \eqref{Side constraint}, \eqref{Integral constraints} and boundary conditions $B_\mathrm{l}(u)=0$, where

\begin{equation}
B_\mathrm{l}\begin{pmatrix} \eta \\ \omega \\ \tau^\prime \\ \zeta^\prime \\ \tau \\ \zeta \end{pmatrix}
=
\begin{pmatrix}
\tau_y^\prime(\tfrac{1}{\beta_0}) \\[1em]
(\mu(1)+\dot{\mu}(1))\tau_y(-\tfrac{1}{\beta_0}) \\[1em]
(\mu(1)+\dot{\mu}(1))\tau_y(0)-\tau^\prime_y(0)\\[0.5em]
\displaystyle \frac{1}{\mu(1)}\zeta(0) - \zeta^\prime(0) + (\mu(1)-1)\omega \\[1em]
\end{pmatrix}.
\label{Explicit formula for Bl}
\end{equation}
An explicit calculation shows that  a complex number $\lambda$ is an eigenvalue of the linear problem (that is, the linear problem
admits a solution of the form $u(x) = \exp(\lambda x) w$ with $w \neq 0$) if and only
if $\lambda=\beta_0 \sigma$, where
\begin{align*}
(\mu(1)-1)^2(\mu(1)&+\dot{\mu}(1))\sigma\tilde{\sigma}\sin\sigma\sin\tilde{\sigma} \\
& =(\sigma^2\beta_0-\alpha_0)\big(\sigma\sin\sigma\cos\tilde{\sigma} + (\mu(1)+\dot{\mu}(1))\tilde{\sigma}\sin\tilde{\sigma}\cos\sigma\big)
\end{align*}
and
$$\tilde{\sigma}=\sigma\sqrt{\frac{\mu(1)}{\mu(1)+\dot{\mu}(1)}}.$$
In particular, $0$ is not an eigenvalue of $L$, and $L$ has a finite number of purely imaginary eigenvalues.

\begin{figure}[h]
\hspace{5cm}\includegraphics[width=6cm]{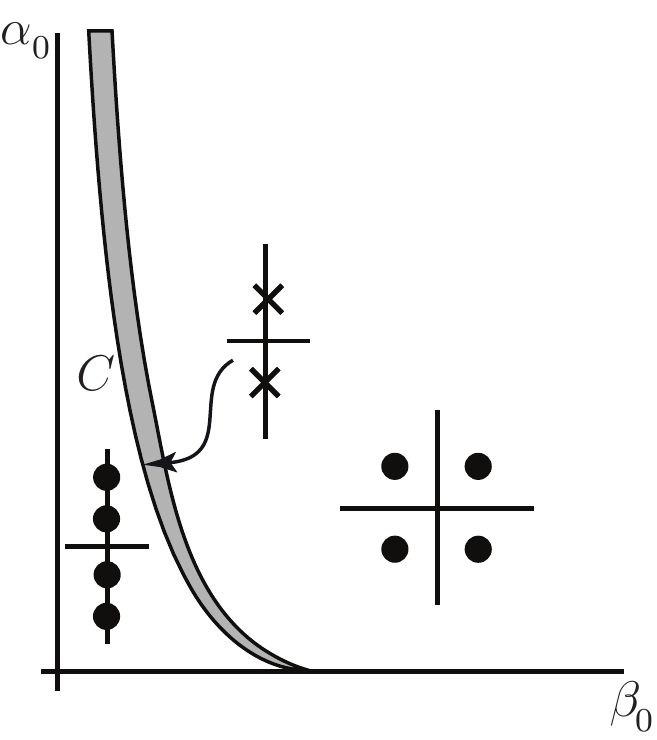}
{\it\caption{Purely imaginary eigenvalues as functions of the parameters $\beta_0$ and $\alpha_0$. A Hamiltonian Hopf bifurcation takes place at points of the curve $C$, and homoclinic solutions exist in the shaded region.} \label{Eigenvalues}}
\end{figure}

According to this calculation, a purely imaginary number $\sigma= \i s$ is an eigenvalue of $L$ if and only if $s=\beta_0 q$, where
$$
q^2\mu(1)(\mu(1)-1)^2  = (\alpha_0+q^2\beta_0)\left( \tilde{q}\coth \tilde{q} + q \mu(1)\coth q\right)
$$
and
$$\tilde{q}=q\sqrt{\frac{\mu(1)}{\mu(1)+\dot{\mu}(1)}}$$
($0$ is not an eigenvalue of $L$). This equation has either zero, one or two pairs $\pm\i q$ of solutions
(see Figure \ref{Eigenvalues}).
A \emph{Hamiltonian Hopf bifurcation} takes place at points of the curve
$$C=\{(\beta_\mathrm{HH}(q),\alpha_\mathrm{HH}(q)): q \in (0,\infty)\}$$
in the $(\beta_0,\alpha_0)$-plane: at the point $(\beta_\mathrm{HH}(q),\alpha_\mathrm{HH}(q))$
two pairs of simple, purely imaginary eigenvalues become complex by
colliding on the imaginary axis at $\pm \i \beta_0 q$ and forming two Jordan chains of
length $2$. Here
\begin{align*}
\beta_\mathrm{HH}(q) & = \frac{\mu(1)(\mu(1)-1)^2}{2(\tilde{q}\coth\tilde{q}+q\mu(1)\coth q)}
+\frac{\mu(1)(\mu(1)-1)^2(\tilde{q}^2\cosech^2\tilde{q}+q^2\mu(1)\cosech^2q)}{2(\tilde{q}\coth\tilde{q}+q\mu(1)\coth q)^2}, \\
\alpha_\mathrm{HH}(q) & = \frac{q^2\mu(1)(\mu(1)-1)^2}{2(\tilde{q}\coth\tilde{q}+q\mu(1)\coth q)}
-\frac{q^2\mu(1)(\mu(1)-1)^2(\tilde{q}^2\cosech^2\tilde{q}+q^2\mu(1)\cosech^2q)}{2(\tilde{q}\coth\tilde{q}+q\mu(1)\coth q)^2}.
\end{align*}
This eigenvalue collision is associated with the bifurcation of a branch of homoclinic solutions
into the region with complex eigenvalues (the shaded region in Figure \ref{Eigenvalues}); we therefore
set $(\beta_0,\alpha_0) =(\beta_\mathrm{HH}(q),\alpha_\mathrm{HH}(q))$ (and $\gamma_0=\alpha_0\beta_0$)
for some $q \in (0,\infty)$.

Choose $e$, $f$ such that 
$$
L e=\i \beta_0 q e,\quad L \overline{e}=-\i \beta_0 q\overline{e}, \qquad
(L-\i \beta_0 q I)f=e,\quad (L+\i \beta_0 q I)\overline{f}=\overline{e},
$$
where $Re = \overline{e}$, $Rf = -\overline{f}$, the `symplectic products' 
$\Omega(e,\overline{f})$ and $\Omega(f,\overline{e})$ are respectively $1$ and $-1$,  and the symplectic products
of all other combinations are zero (note that $\Omega$ acts bilinearly on pairs of complex vectors).

\begin{theorem}
The Hamiltonian formulation of our ferrohydrostatic problem admits a locally invariant manifold
of the form $\{A e + B f + \overline{A} \overline{e} + \overline{B} \overline{f} + \tilde{r}(A,B,\overline{A},\overline{B},\varepsilon)\}$, where
$\tilde{r}(A,B,\overline{A},\overline{B},\varepsilon)$ satisfies
$\tilde{r} = O(|(A,B,\overline{A},\overline{B})||(\varepsilon,A,B,\overline{A},\overline{B})|)$. The flow on this manifold is described by the
two-degree-of-freedom Hamiltonian system
\begin{equation}
A_x = \frac{\partial \tilde{H}_\mathrm{f}^\varepsilon}{\partial \overline{B}}, \qquad B_x = -\frac{\partial \tilde{H}_\mathrm{f}^\varepsilon}{\partial \overline{A}}, \label{Reduced system}
\end{equation}
where
$$\tilde{H}_\mathrm{f}^\varepsilon(A,B,\overline{A},\overline{B}) = H_\mathrm{f}^\varepsilon(A e + B f + \overline{A} \overline{e} + \overline{B} \overline{f} + \tilde{r}(A,B,\tilde{A},\overline{B},\varepsilon)),$$
which is reversible with reverser $R(A,B) =(\overline{A},-\overline{B})$. In particular,
\begin{itemize}
\item[(i)]
all small bounded solutions of \eqref{HE}, \eqref{NLBC} (with constraints \eqref{Side constraint}, \eqref{Integral constraints}) take the form
\begin{equation}
u(x)=A(x) e + B(x) f + \overline{A}(x) \overline{e} + \overline{B}(x) \overline{f} + \tilde{r}(A(x),B(x),\tilde{A}(x),\tilde{B}(x),\varepsilon), \label{Full vs reduced}
\end{equation}
where $(A(x),B(x))$, $x \in {\mathbb R}$ solves \eqref{Reduced system};
\item[(ii)]
any small, bounded solution $(A(x),B(x))$, $x \in {\mathbb R}$ of \eqref{Reduced system}
defines a small, bounded solution \eqref{HE}, \eqref{NLBC} (with constraints \eqref{Side constraint}, \eqref{Integral constraints})
via formula \eqref{Full vs reduced}.
\end{itemize}
\end{theorem}

\subsection{Proof of the reduction theorem} \label{Reduction theory}

The key result is the following theorem, which is 
a parametrised,
Hamiltonian version of a reduction principle for quasilinear 
evolutionary equations
presented by Mielke \cite[Theorem 4.1]{Mielke88a} (see Buffoni,
Groves \& Toland \cite[Theorem 4.1]{BuffoniGrovesToland96}).

\begin{theorem}\label{Mielke's theorem}
Consider the differential equation
\begin{equation}
u_x = \LL u + \NN(u;\lambda), \label{feq}
\end{equation}
which represents Hamilton's equations for the reversible Hamiltonian system 
$(M,\Omega^\lambda,H^\lambda)$.
Here $u$ belongs to a Hilbert space $\XX$, $\lambda \in {\mathbb R}^\ell$ is
a parameter and $\LL: \DD(\LL) \subset \XX \rightarrow \XX$ is a densely 
defined, closed
linear operator. Regarding $\DD(\LL)$ as a Hilbert space equipped with 
the graph norm,
suppose that $0$ is an equilibrium point of \eqref{feq} when $\lambda=0$ and that
\begin{enumerate}
\item[(H1)] The part of the spectrum $\sigma(\LL)$ of $\LL$ which lies on 
the imaginary
axis consists of a finite number of eigenvalues of finite multiplicity
and is separated from the
rest of $\sigma(\LL)$ in the sense of Kato, so that $\XX$ admits the
decomposition $\XX=\XX_1\oplus \XX_2$, where $\XX_1=\PP(\XX)$, 
$\XX_2=(I-\PP)(\XX)$ and
$\PP$ is the spectral projection corresponding the purely imaginary part
of $\sigma(\LL)$.
\item[(H2)] The operator $\LL_2=\LL|_{\XX_2}$ satisfies the estimate
$$\|(\LL_2-\i s I)^{-1}\|_{\XX_2 \rightarrow \XX_2} \leq \frac{C}{1+|s|}, \qquad s \in {\mathbb R},$$
for some constant $C$ that is independent of $s$.
\item[(H3)] There exists a natural number $k$ and neighbourhoods 
$\Lambda \subset
{\mathbb R}^\ell$ of $0$ and $U \subset \DD(\LL)$ of $0$ such that $\NN$ is $(k+1)$
times continuously differentiable on $U \times \Lambda$, its 
derivatives are bounded
and uniformly continuous on $U \times \Lambda$ and $\NN(0,0)=0$,
$\mathrm{d}_1\NN[0,0]=0$.
\end{enumerate}
Under these hypotheses there exist neighbourhoods $\tilde{\Lambda} \subset
\Lambda$ of $0$ and $\tilde{U}_1 \subset U \cap \XX_1$, $\tilde{U}_2 
\subset U \cap \XX_2$
of $0$ and a reduction function
$r:\tilde{U}_1\times\tilde\Lambda\to \tilde{U}_2$ with the following properties.
The reduction function $r$ is $k$ times continuously differentiable
on $\tilde{U}_1\times\tilde\Lambda$, its derivatives are bounded and
uniformly continuous on $\tilde{U}_1\times\tilde\Lambda$ and
$r(0;0)=0$, $\mathrm{d}_1r[0;0]=0$. The graph
$\tilde{M}^\lambda=\{u_1+r(u_1;\lambda) \in \XX_1 \oplus \XX_2: u_1 \in \tilde{U}_1\}$
is a Hamiltonian centre manifold for \eqref{feq}, so that
\begin{list}{(\roman{count})}{\usecounter{count}}
\item
$\tilde{M}^\lambda$ is a locally invariant manifold of \eqref{feq}:
through every point in $\tilde{M}^\lambda$  there passes a unique solution
of \eqref{feq} that remains on $\tilde{M}^\lambda$ as long as it
remains in $\tilde{U}_1\times \tilde{U}_2$.
\item
Every small bounded solution $u(x)$, $x\in{\mathbb R}$ of
\eqref{feq} that satisfies $(u_1(x),u_2(x))\in \tilde{U}_1\times \tilde{U}_2$
lies completely in $\tilde{M}^\lambda$.
\item
Every solution $u_1:(x_1,x_2) \to \tilde{U}_1$ of the reduced equation
\begin{equation}
u_{1x}=\LL u_1+\PP\NN(u_1+r(u_1;\lambda);\lambda) \label{req}
\end{equation}
generates a solution
\begin{equation}
u(x)=u_1(x)+r(u_1(x);\lambda) \label{waves}
\end{equation}
of the full equation \eqref{feq}.
\item
$\tilde{M}^\lambda$ is a symplectic submanifold of $M$ and the flow determined
by the Hamiltonian system $(\tilde{M}^\lambda,\tilde{\Omega}^\lambda,
\tilde{H}^\lambda)$,
where the tilde denotes  restriction to $\tilde{M}^\lambda$, coincides
with the flow on $\tilde{M}^\lambda$ determined by
$(M,\Omega^{\lambda},H^\lambda)$.
The reduced equation \eqref{req} is reversible and represents Hamilton's equations for
$(\tilde{M}^\lambda,\tilde{\Omega}^\lambda, \tilde{H}^\lambda)$.
\end{list}
\end{theorem}

Mielke's theorem cannot be applied directly to the Hamiltonian system
$(M_0,\Omega_\mathrm{f},H_\mathrm{f}^\varepsilon)$ because of the
nonlinear boundary conditions \eqref{NLBC} in the domain of the Hamiltonian vector field $v_{H^\varepsilon_\mathrm{f}}$
(the right-hand sides of \eqref{HE 1}--\eqref{HE 6} define a smooth mapping
$f^\varepsilon: N_1 \rightarrow M_0$ with
$v_{H^\varepsilon_\mathrm{f}}|_u=f^\varepsilon(u)$ for any $u \in \DD(v_{H^\varepsilon_\mathrm{f}})$).
We overcome this difficulty by using the change of variable
$(\eta,\omega,\tau^\prime,\zeta^\prime,\tau,\zeta)=G(\eta,\omega,\chi^\prime,\xi^\prime,\chi,\xi)$, where
\begin{align*}
\tau^\prime & = \chi^\prime + \int_0^y \left\{\frac{K_1^\prime W\xi^\prime}{\sqrt{1-W^2}}+K_2^\prime\chi_y^\prime - \chi_y^\prime \right\}\dt - \beta_0\hat{\tau}, \\
\zeta^\prime & = \xi^\prime - \frac{WK_1^\prime\chi_y^\prime}{\sqrt{1-W^2}}+K_2^\prime\xi^\prime
-(\mu(1)-1)\frac{W}{\sqrt{1-W^2}} - \xi^\prime + (\mu(1)-1)\omega - 3y^2\beta_0^3 \hat{\zeta}, \\
\tau & = \chi + \frac{1}{\mu(1)+\dot{\mu}(1)}\!\int_y^0\!\! \left\{\!-\frac{K_1W\xi}{\sqrt{1-W^2}}-\mu^\dagger(K_2\chi_y+1)
+\mu(1)+(\mu(1)+\dot{\mu}(1))\chi_y \right\}\! \dt - \beta_0\hat{\tau}, \\
\zeta & = \xi + \mu(1) \left(-\frac{WK_1\chi_y}{\sqrt{1-W^2}}+\frac{K_2}{\mu^\dagger}\xi - \frac{1}{\mu(1)}\xi\right) - 3y^2\beta_0^3 \hat{\zeta},
\end{align*}
and
{\allowdisplaybreaks\begin{align*}
\hat{\tau} & =
\frac{1}{2(\mu(1)+\dot{\mu}(1))}
\!\int_{-\tfrac{1}{\beta_0}}^0 \int_y^0 \!\!\left\{\!-\frac{K_1W\xi}{\sqrt{1-W^2}}-\mu^\dagger(K_2\chi_y+1)
+\mu(1)+(\mu(1)+\dot{\mu}(1))\chi_y\! \right\}\!\dt\dy \\
& \qquad\mbox{}+
\frac{1}{2}\int_0^{\tfrac{1}{\beta_0}} \int_0^y\! \left\{\frac{K_1^\prime W\xi^\prime}{\sqrt{1-W^2}}+K_2^\prime\chi_y^\prime - \chi_y^\prime \right\}\dt\dy, \\
\hat{\zeta} & = \frac{\mu(1)}{2}\int_0^{\tfrac{1}{\beta_0}} \left\{-\frac{WK_1\chi_y}{\sqrt{1-W^2}}+\frac{K_2}{\mu^\dagger}\xi - \frac{1}{\mu(1)}\xi\right\}\dy \\
& \qquad\mbox{}+\frac{1}{2}\int_0^{\tfrac{1}{\beta_0}}\left\{
-\frac{WK_1^\prime\chi_y^\prime}{\sqrt{1-W^2}}+K_2^\prime\xi^\prime
-(\mu(1)-1)\frac{W}{\sqrt{1-W^2}} - \xi^\prime + (\mu(1)-1)\omega\right\}\dy,
\end{align*}}
which transforms the nonlinear boundary conditions in $\DD(v_{H^\varepsilon_\mathrm{f}})$ into their linearisations
$$B_\mathrm{l}(\eta,\omega,\tau^\prime,\zeta^\prime,\tau,\zeta)=0$$
where $B_\mathrm{l}=\mathrm{d}B[0]$ (see formula \eqref{Explicit formula for Bl}).

\begin{lemma}$ $
\begin{itemize}
\item[(i)]
There exists an open neighbourhood $\hat{N}_1$ of the origin such that
$G: N_1 \rightarrow \hat{N}_1$ is a diffeomorphism.
\item[(ii)]
For each $u \in N_1$ the operator $\mathrm{d}G[u]: M_1 \rightarrow M_1$
extends to an isomorphism $\widetilde{\mathrm{d}G}[u]: M_0 \rightarrow M_0$.
The operators $\widetilde{\mathrm{d}G}[u]$, $\widetilde{\mathrm{d}G}[u]^{-1} \in {\mathcal L}(M_0,M_0)$
depend smoothly upon $u \in N_1$.
\end{itemize}
\end{lemma}
{\bf Proof.} (i) We observe that $G$ is a smooth, near identity mapping $N_1 \rightarrow M_1$ and apply the
inverse-function theorem. (ii) The result for $\mathrm{d}G[u]$ follows from a direct calculation, while the second assertion
is proved using the argument given by Groves \& Mielke \cite[Lemma 3.3]{GrovesMielke01}.\qed\medskip

A simple calculation shows that the diffeomorphism $G$ transforms
$$u_x = f^\varepsilon(u)$$
into
\begin{equation}
v_x=\hat{f}^\varepsilon(v), \label{HE with linear BC}
\end{equation}
where $\hat{f}^\varepsilon:\hat{N}_1 \to M_0$ is the smooth vector field defined by
$$
\hat{f}^\varepsilon(v)=\widetilde{{\mathrm d}G}\,[G^{-1}(v)]\,
(f^\varepsilon(G^{-1}(v))).
$$
Formula \eqref{HE with linear BC} represents Hamilton's equations for the Hamiltonian
system $(M_0,\Gamma,\hat{H}^\varepsilon_\mathrm{f})$, where
$$
\Gamma\big|_{n}(v_1,v_2)=
\Omega(\widetilde{{\mathrm d}G}\,[G^{-1}(v)]^{-1}(v_1),
\widetilde{{\mathrm d}G}\,[G^{-1}(n)]^{-1}(v_2)), \qquad n \in \hat{N}_1, \ v_1,v_2 \in TX_0|_v,
$$
and
$$
\hat{H}^\varepsilon_\mathrm{f}(n)=H^\varepsilon_\mathrm{f}(G^{-1}(n)), \qquad n \in \hat{N}_1.
$$
The domain of the Hamiltonian vector field $v_{\hat{H}^\varepsilon_\mathrm{f}}$ is
$$\DD(v_{\hat{H}^\varepsilon_\mathrm{f}})=\{(\eta,\omega,\tau^\prime,\zeta^\prime,\tau,\zeta) \in \hat{N}_1: B_\mathrm{l}(\eta,\omega,\tau^\prime,\zeta^\prime,\tau,\zeta)=0\}$$
and
$v_{\hat{H}^\varepsilon_\mathrm{f}}|_{n}=\hat{f}^\varepsilon(n)$ for any $n\in \DD(v_{\hat{H}^\varepsilon_\mathrm{f}})$.

The next step is to verify that \eqref{HE with linear BC} satisfies the hypotheses
of Theorem \ref{Mielke's theorem} (with $\XX=M_0$), so that we obtain
a finite-dimensional reduced Hamiltonian system
$(\tilde{M}^\varepsilon,\tilde{\Gamma}^\varepsilon,\tilde{H}_\mathrm{f}^\varepsilon)$.
We write \eqref{HE with linear BC} as
$$
u_{x} = L u + N^\varepsilon(u),
$$
in which the linear operator $L: \DD(L) \subset M_0 \rightarrow M_0$ with
$L={\mathrm d} g^0[0]$ and
$$\DD(L)=\{(\eta,\omega,\chi^\prime,\xi^\prime,\chi,\xi) \in N_1: B_\mathrm{l}(\eta,\omega,\chi^\prime,\xi^\prime,\chi,\xi)=0\}$$
is given explicitly by the formula \eqref{Explicit formula for L}.
(Observe that ${\mathrm d} g^0[0]={\mathrm d} f^0[0]$ and $B_\mathrm{l}=\mathrm{d}B[0]$, so that
$L$ is the formal linearisation of $v_{H^\varepsilon_\mathrm{f}}$.)
It follows from Lemma \ref{lin} below that $L$ satisfies hypotheses (H1) and (H2);
hypothesis (H3) is clearly satisfied for an arbitrary value of $k$. Part (i) of Lemma \ref{lin}
is proved using the elementary theory of ordinary differential equations, while part
(ii) established using arguments similar to those employed for other problems
treated using centre-manifold reduction (e.g.\ see Buffoni, Groves \& Toland
\cite[Proposition 3.2]{BuffoniGrovesToland96} or Groves \& Wahl\'{e}n \cite[Lemma 3.4]{GrovesWahlen07}).

\begin{lemma}\label{lin}
\hfill
\begin{list}{(\roman{count})}{\usecounter{count}}
\item
The spectrum $\sigma(L)$ of $L$ consists entirely of isolated
eigenvalues of finite algebraic multiplicity.
A complex number $\lambda$ is an eigenvalue of $L$ if and only if $\lambda=\beta_0 \sigma$, where
\begin{align*}
(\mu(1)-1)^2(\mu(1)&+\dot{\mu}(1))\sigma\tilde{\sigma}\sin\sigma\sin\tilde{\sigma} \\
& =(\sigma^2\beta_0-\alpha_0)\big(\sigma\sin\sigma\cos\tilde{\sigma} + (\mu(1)+\dot{\mu}(1))\tilde{\sigma}\sin\tilde{\sigma}\cos\sigma\big)
\end{align*}
and
$$\tilde{\sigma}=\sigma\sqrt{\frac{\mu(1)}{\mu(1)+\dot{\mu}(1)}}.$$
(In particular, $0 \not\in \sigma(L)$ and $\sigma(L) \cap \i {\mathbb R}$ is a finite set.)
\item
There exist real constants $C$, $s_0>0$ such that
$$
\|(L-\i sI)^{-1}\|_{\LL(M_0,M_0)}\leq\frac{C}{|s|}
$$
for each real number $s$ with $|s|>s_0$.
\end{list}
\end{lemma}

The centre manifold  $\tilde{M}^\varepsilon$
is equipped with the single coordinate chart $\tilde{U}_1 \subset \XX_1$ and  coordinate
map $\pi: \tilde{M}^\varepsilon \rightarrow \tilde{U}_1$ defined by
$\pi^{-1}(u_1) = u_1 + r(u_1;\varepsilon)$. It is however more convenient to use
an alternative coordinate map for calculations.
According to the parameter-dependent version of Darboux's theorem (e.g.\ see
Buffoni \& Groves \cite[Theorem 4]{BuffoniGroves99})
there exists a near-identity change of variable
$$
u_1 = \hat{u}_1 + D(\hat{u}_1;\varepsilon)
$$
of class $C^{k-1}$ which transforms $\tilde{\Gamma}^\varepsilon$
into $\Upsilon$, where
$$\Upsilon(v_1,v_2) = \tilde{\Gamma}^0|_0(v_1,v_2)=\Gamma|_0(v_1,v_2).$$
Define the function
$\tilde{r}: U_1 \times \tilde{\Lambda} \rightarrow \tilde{U}_1 \times \tilde{U}_2$
with $U_1=\PP G^{-1}(\tilde{U}_1\times \tilde{U}_2)$
(which in general has components in $\XX_1$ and $\XX_2$) by the formula
\begin{equation}
\hat{u}_1 + \tilde{r}(\hat{u}_1;\varepsilon) = G^{-1}\big(\hat{u}_1+D(\hat{u}_1;\varepsilon)
+r(\hat{u}_1+D(\hat{u}_1;\varepsilon);\varepsilon)\big), \label{Definition of tilde r}
\end{equation}
where $\tilde{r}(0;0)=0$, $\mathrm{d}_1\tilde{r}[0;0]=0$, and equip $\tilde{M}^\varepsilon$
with the coordinate map $\hat{\pi}: \tilde{M}^\varepsilon \rightarrow U_1$ given by
$\hat{\pi}^{-1}(\hat{u}_1) = \hat{u}_1 + \tilde{r}(\hat{u}_1;\varepsilon)$.

We use a basis for $\XX_1$ with respect to which $\Upsilon$ is the canonical symplectic $2$-form
(a `symplectic basis'). Choose $e$, $f \in \DD(L)$ such that 
$$
L e=\i \beta_0 q e,\quad L \overline{e}=-\i \beta_0 q\overline{e}, \qquad
(L-\i \beta_0 q I)f=e,\quad (L+\i \beta_0 q I)\overline{f}=\overline{e},
$$
where $Re = \overline{e}$, $Rf = -\overline{f}$ and
$\Omega(e,\overline{f})=1$, $\Omega(f,\overline{e})=-1$ and the symplectic products
of all other combinations are zero (note that $\Omega$ acts bilinearly on pairs of complex vectors).
It follows that $\{e,f,\overline{e},\overline{f}\}$ is a symplectic basis for the central subspace of $L$
(so that
the coordinates $A$, $B$, $\overline{A}$ and $\overline{B}$ in the $e$, $f$, $\overline{e}$ and $\overline{f}$ directions
are canonical coordinates) and the action of the reverser $R$ on this
space is given by
$$R(A,B) =(\overline{A},-\overline{B}).$$
We can now identify $(\tilde{M}^\varepsilon,\tilde{\Gamma}^\varepsilon,\tilde{H}_\mathrm{f}^\varepsilon)$ with
the four-dimensional canonical Hamiltonian
system $(U_1, \Upsilon, \tilde{H}_\mathrm{f}^\varepsilon)$ using the coordinate map $\hat{\pi}$, where
$$
\Upsilon((A^1,B^1,\overline{A^1},\overline{B^1}),
(A^2,B^2,\overline{A^2},\overline{B^2}))
=A^1\overline{B^2}-A^2\overline{B^1}
+\overline{A^1}B^2-\overline{A^2}B^1
$$
and
$$\tilde{H}^\varepsilon_\mathrm{f}(A,B) = \hat{H}^\varepsilon_\mathrm{f}(\hat{u}_1+\tilde{r}(\hat{u}_1;\varepsilon)), \qquad
\hat{u}_1 = A e + B f + \overline{A}\overline{e} + \overline{B}\overline{f};
$$
Hamilton's equations for $(U_1, \Upsilon, \tilde{H}_\mathrm{f}^\varepsilon)$ are given by \eqref{Reduced system}.

\section{Homoclinic bifurcation} \label{NF theory}

The coordinates $A$, $B$, $\overline{A}$, $\overline{B}$ for the reduced Hamiltonian system may be chosen so that
the reduced Hamiltonian takes the normal form
\begin{eqnarray}
\lefteqn{\tilde{H}_\mathrm{f}^\varepsilon(A,B)=\i \beta_0 q(A\overline{B}-\overline{A}B)+|B|^2}\hspace{0.4in} \nonumber\\
& & \qquad\mbox{} +H_\mathrm{NF}(|A|^2,\i(A\overline{B}-\overline{A}B),\varepsilon)
+ O(|(A,B)|^2|(\varepsilon,A,B)|^{n_0}), \label{NF Hamiltonian}
\end{eqnarray}
where $H_\mathrm{NF}$ is a real polynomial of order $n_0+1$
satisfying
$$H_\mathrm{NF}(|A|^2,\i(A\overline{B}-\overline{A}B),\varepsilon)=O(|(A,B)|^2|(\varepsilon,A,B)|)$$
and $n_0 \geq 2$ is a fixed integer; Hamilton's equations for the reduced system
are given by
\begin{eqnarray}
A_x& = &
\i \beta_0 q A+B+ \i A\partial_2 H_\mathrm{NF}(|A|^2,\i
(A\overline B-\overline AB),\varepsilon)+O(|(A,B)||(\varepsilon,A,B)|^{n_0}), \label{IP BG 1} \\
B_x& = &
\i \beta_0 q B+\i B\partial_2 H_\mathrm{NF}(|A|^2,\i
(A\overline B-\overline AB),\varepsilon) \nonumber \\
& & 
\quad\ \mbox{}-A\partial_1 H_\mathrm{NF}(|A|^2,\i
(A\overline B-\overline AB),\varepsilon)+O(|(A,B)||(\varepsilon,A,B)|^{n_0}). \label{IP BG 2} 
\end{eqnarray}

Existence theories for homoclinic solutions to \eqref{IP BG 1}, \eqref{IP BG 2} have been given by
Iooss \&\linebreak
P\'{e}rou\`{e}me \cite{IoossPeroueme93} and Buffoni \& Groves \cite{BuffoniGroves99}
under the assumption that the coefficients $c_1$ and $c_3$ in the expansion
\begin{eqnarray*}
\lefteqn{H_\mathrm{NF} = \varepsilon c_1 |A|^2 + \varepsilon \i c_2(A\overline{B}-\overline{A}B)
+c_3|A|^4} \hspace{0.5in} \\
& & \mbox{}+\i c_4|A|^2(A\overline{B}-\overline{A}B)-c_5(A\overline{B}-\overline{A}B)^2
+ \varepsilon^2 c_6 |A|^2 + \varepsilon^2\i c_7 (A\overline{B}-\overline{A}B)+\ldots
\end{eqnarray*}
are respectively negative and positive.

\begin{theorem}
Suppose that $c_1<0$ and $c_3>0$.
\begin{list}{(\roman{count})}{\usecounter{count}}
\item (Iooss \& P\'{e}rou\`{e}me) For each sufficiently small,
positive value of $\varepsilon$ the 
two-degree-of-freedom Hamiltonian system \eqref{IP BG 1}, \eqref{IP BG 2}
has two distinct symmetric homoclinic 
solutions.
\item (Buffoni \& Groves) For each sufficiently small, positive value of $\varepsilon$ the 
two-degree-of-freedom Hamiltonian system \eqref{IP BG 1}, \eqref{IP BG 2}
has an infinite number of geometrically 
distinct homoclinic solutions which generically resemble multiple
copies of one of the homoclinic solutions in part (i).
\end{list}
The homoclinic solutions identified above correspond to envelope patterns whose amplitude
is $O((-c_1\varepsilon)^{1/2})$ and which decay exponentially as $x \rightarrow \pm\infty$;
they are sketched in Figure \ref{Interface sketches}. 
\end{theorem}

The coefficients $c_1$ und $c_3$ are given by the formulae
$$c_1 = 2H_2^1[e,\overline{e}]$$
and
$$c_3=6H_4^0[e,e,\overline{e},\overline{e}]+3H_3^0[\tilde{r}^0_{2000},\overline{e},\overline{e}]+3H_3^0[\tilde{r}^0_{1010},e,\overline{e}],$$
where
$\tilde{r}^j_{k_1k_2k_3k_4}$ is the coefficient of
$\varepsilon^jA^{k_1}B^{k_2}\overline{A}^{k_3}\overline{B}^{k_4}$ in the
Taylor expansion of $\tilde{r}$ and
$$H_k^0 = \frac{1}{k!}\mathrm{d}^kH^0[0], \qquad
H_k^1 = \frac{1}{k!}\mathrm{d}^k \left(\left.\frac{\partial H^\varepsilon}{\partial \varepsilon}\right|_{\varepsilon=0}\right)[0];
$$
the quantities $\tilde{r}^0_{2000}$ and $\tilde{r}^0_{1010}$
are found from the linear boundary-value problems
\begin{eqnarray*}
(f_1^0-2\i \beta_0 qI)\tilde{r}^{0}_{2000} & = & -f_2^0[e,e], \\
B_\mathrm{l} \tilde{r}^{0}_{2000} & = & -B_2[e,e], \\
\\
f_1^0\tilde{r}^{0,0}_{1010} & = & -2f_2^0[e,\overline{e}], \\
B_\mathrm{l} \tilde{r}^{0}_{1010} & = & -2B_2[e,\overline{e}],
\end{eqnarray*}
where
$$f_k^0 = \frac{1}{k!}\mathrm{d}^kf^0[0], \qquad B_k = \frac{1}{k!}\mathrm{d}^kB[0].$$
These formulae are derived using the method explained by Groves \& Mielke \cite[Appendix B]{GrovesMielke01}.

Attempting to compute explicit general expressions for $c_1$ and $c_2$ leads to unwieldy formulae (it appears more
appropriate to calculate them numerically for a specific choice of $\mu$, that is a specific magnetisation law).
Here we confine ourselves to stating the values of the coefficients for two particular special cases.\newline

\begin{figure}[h]
\centering
\includegraphics[scale=0.4]{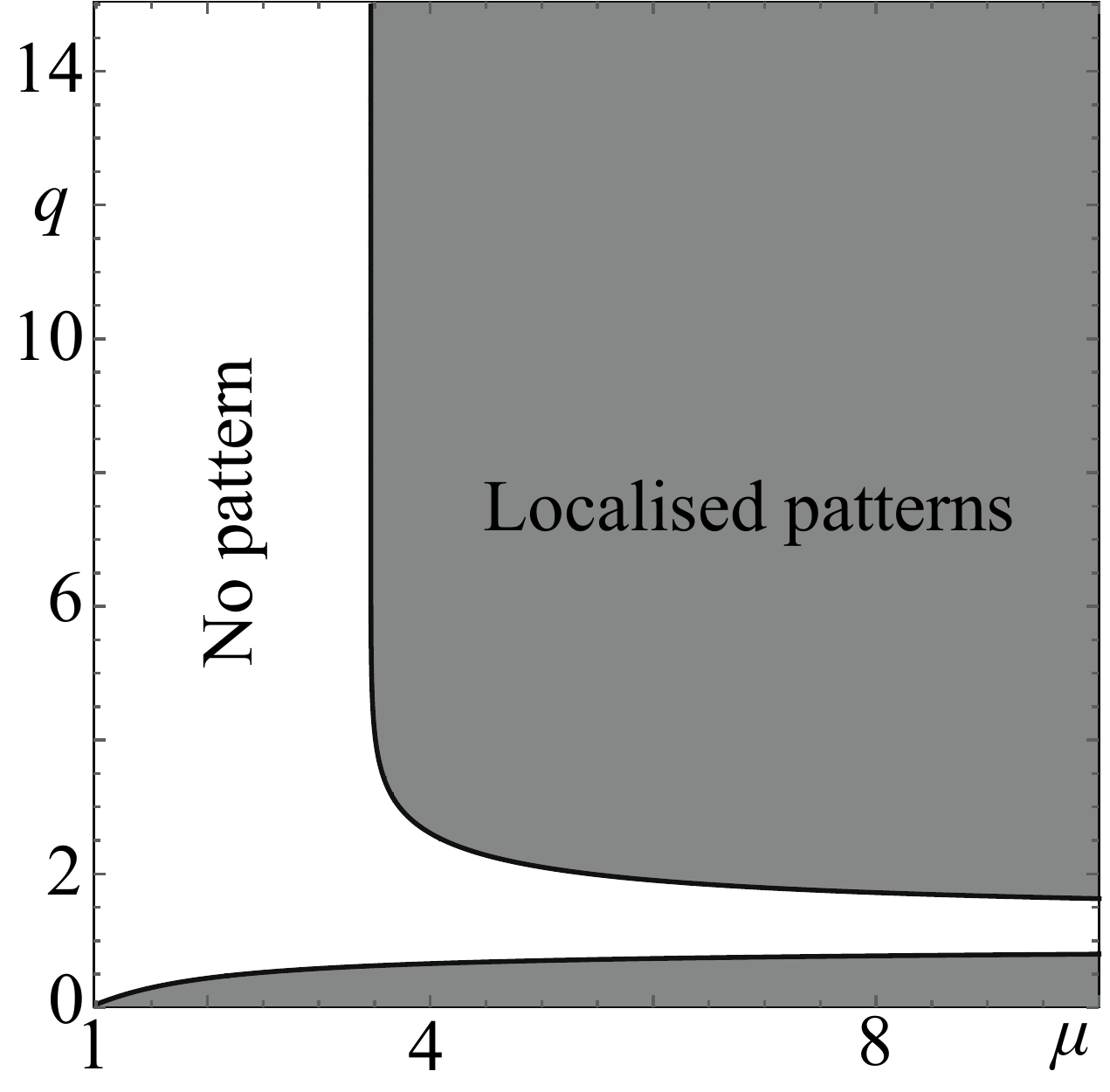}\hspace{2cm}\includegraphics[scale=0.4]{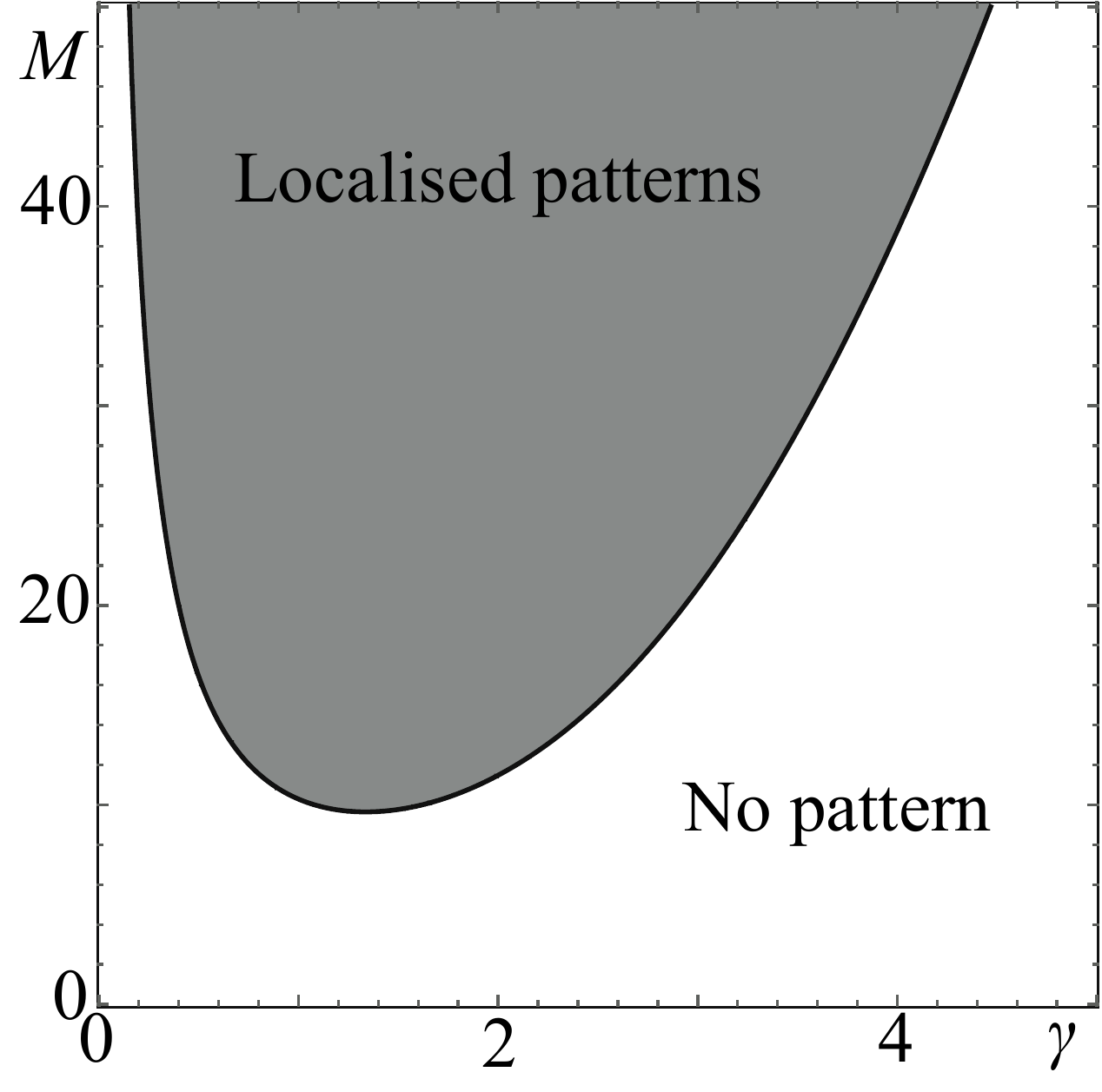}
{\it \caption{The sign of the coefficient $c_3$ as a function of $\mu$ and $q$ for a linear magnetisation law (left) and as a function of
$M$ and $\gamma$ for the Langevin magnetisation law \eqref{Langevin} in a ferrofluid of great depth (right). The coefficient $c_1$ is negative in both panels. The shaded areas show the regions in which localised one-dimensional interfaces bifurcate from the trivial state.}
\label{c3 sign}}
\end{figure}

\noindent\emph{(i) Constant relative permeability $\mu$ (corresponding to a linear magnetisation law)}: We find that
\small
\begin{align*}
c_1 = & -\left(1 + \frac{\mu}{\beta_0}\frac{(\mu-1)^2}{\mu+1} (q \tanh q-1)\sech^2 q  \right)^{\!\!-1}, \\[2mm]
c_3 = &\frac{1}{2}\left(1 + \frac{\mu}{\beta_0}\frac{(\mu-1)^2}{\mu+1} (q \tanh q-1)\sech^2 q \right)^{\!\!-1} \\
& \times\!\! \Bigg(\frac{q^4\beta_0^4\mu^2}{4}\frac{(\mu-1)^6}{(\mu+1)^4} \frac{(-3-4\cosh 2q + \cosh 4q)(-2+\sech^2 q + 4 \sech 2q)}
{\left(2q\beta_0\sdfrac{(\mu-1)^2}{\mu+1}\mu \tanh 2q- \gamma_0 - 4q^2\beta_0^2\right)\cosh^2 q \cosh 2q} \\
& \qquad\!\mbox{}+\frac{q^4\beta_0^4 \mu^2}{\gamma_0}\frac{(\mu-1)^6}{(\mu+1)^4} \sech^4 q \\
& \qquad\!\mbox{}+\frac{q^3\beta_0^3}{4}\frac{\cosech 2q}{(\mu+1)^3}\! \left(\!-16\mu(\mu^2-1)^2\!\cosh^2 q + 16(\mu-1)^4\mu\sech 2q+6q \beta_0 (\mu+1)^3\!\sinh 2q \right)\!\!\!\Bigg),
\end{align*}
\normalsize
\noindent where $\beta_0 q$ is the critical wavenumber associated with the Rosensweig instability.
The sign of $c_3$ clearly depends upon $\mu$ and $q$ (see Figure \ref{c3 sign}). For large fluid depths 
one requires $\mu\geq\mu_{\mathrm{c}}\approx 3.5$ for localised one-dimensional free surfaces to bifurcate from the trivial state, while for experimentally relevant ferrofluids ($2<\mu<6$) one requires $q>2$. For sufficiently shallow depths, we find that bifurcation occurs for all values of $\mu$. However, due to the choice of the boundary conditions, the model described in Section \ref{ferroeqns} may become invalid at these shallow depths.

\noindent\emph{(ii) Small values of $\beta_0$ (corresponding to deep fluids)}: Abbreviating $\mu(1)$, $\dot{\mu}(1)$,
$\ddot{\mu}(1)$, $\dddot{\mu}(1)$ to respectively $\mu_1$, $\dot{\mu}_1$, $\ddot{\mu}_1$,
$\dddot{\mu}_1$, one finds that
\begin{align*}
c_1 = & -1 +o(1), \\[2mm]
c_3=& \frac{3 s^4}{4}-\frac{4s^5}{(\mu_1-1)^2} \left(1+\frac{m}{\mu_1+\dot{\mu}_1}\right)\\
& \quad\mbox{}-
\frac{s^6}{(\mu_1-1)^4\mu_1(\mu_1+\dot{\mu}_1)^3}\Big(-6 \mu_1^2 \ddot{\mu}_1+6 \mu_1 \ddot{\mu}_1+24 \mu_1^3 \dot{\mu}_1+24 \mu_1^2 \dot{\mu}_1^2-50 \mu_1^2 \dot{\mu}_1+8 \mu_1 \dot{\mu}_1^3\\
&\hspace{2.2in}\mbox{}
-22 \mu_1 \dot{\mu}_1^2+26 \mu_1 \dot{\mu}_1+6 \dot{\mu}_1^2+8 \mu_1^4-16 \mu_1^3+8 \mu_1^2\Big)\\
& \quad\mbox{}+
\frac{ms^7}{4(\mu_1-1)^4\mu_1^3(\mu_1+\dot{\mu}_1)^4}\Big(
4 \mu_1^3 \dddot{\mu}_1-9 \mu_1^2 \ddot{\mu}_1^2+20 \mu_1^3 \ddot{\mu}_1-5 \dot{\mu}_1^4-18 \mu_1 \dot{\mu}_1^3-37 \mu_1^2 \dot{\mu}_1^2\\
&\hspace{2.2in}\mbox{}
+12 \mu_1^3 \dot{\mu}_1+4 \mu_1^2 \dot{\mu}_1 \dddot{\mu}_1-2 \mu_1 \dot{\mu}_1^2 \ddot{\mu}_1-18 \mu_1^2 \dot{\mu}_1 \ddot{\mu}_1\Big) \\
& \quad\mbox{}-\frac{4 s^4 \left(\dot{\mu}_1+\mu_1-1\right)a_1}{(\mu_1-1)^3 \left(\mu_1+\dot{\mu}_1\right)}\\
& \quad\mbox{}+
\frac{s^4a_2}{(\mu_1-1)^3\mu_1(\mu_1+\dot{\mu}_1)^2}\Big(\mu_1^2 \ddot{\mu}_1-\mu_1 \ddot{\mu}_1+11 \mu_1^2 \dot{\mu}_1
+5 \mu_1 \dot{\mu}_1^2-7 \mu_1 \dot{\mu}_1-\dot{\mu}_1^2+4 \mu_1^3-4 \mu_1^2\Big)\\
& \quad\mbox{}+o(1),
\end{align*}
as $\beta_0 \rightarrow 0$, where
$$s=\frac{(\mu_1-1)^2\mu_1}{2(m+\mu_1)}, \qquad m = \sqrt{\frac{\mu_1}{\mu_1+\dot{\mu}_1}}$$
and
\small
\begin{align*}
\begin{pmatrix} a_1 \\ a_2 \end{pmatrix}
=&\begin{pmatrix}2s & & 2ms(\mu_1+\dot{\mu}_1) \\
\\
2 \mu_1 s-\sdfrac{5 s^2}{\mu_1-1} & & 2ms(\mu_1+\dot{\mu}_1) +\sdfrac{5 s^2}{\mu_1-1}
\end{pmatrix}^{\!\!-1} \\[4mm]
&\quad\left(\begin{array}{l}
-\sdfrac{4 (m+1) s^3}{\mu_1-1}-\sdfrac{s^4 \left(\dot{\mu}_1^2+\mu_1 \left(\ddot{\mu}_1+3 \dot{\mu}_1\right)\right)}{(\mu_1-1)^2 \mu_1 \left(\mu_1+\dot{\mu}_1\right)^2} \\
\\
-2s^2(\mu_1-1)\mu_1 - \sdfrac{s^4}{(\mu_1-1)^2\mu_1(\mu_1+\dot{\mu}_1)^2}
\\
\hspace{1.33in}\times\Big(\mu_1 \ddot{\mu}_1-28 \mu_1^2 \dot{\mu}_1 -14 \mu_1 \dot{\mu}_1^2+17 \mu_1 \dot{\mu}_1+\dot{\mu}_1^2-14 \mu_1^3+14 \mu_1^2\Big)
\end{array}\right)\!\!.
\end{align*}
\normalsize
(Note that the critical wavenumber associated with the Rosensweig instability is $s+o(1)$ as $\beta_0 \rightarrow 0$.)
Figure \ref{c3 sign} shows the sign of $c_3$ for the Langevin magnetisation law
\begin{equation}
\mu(s) = 1+\frac{M}{s}\left(\coth(\gamma s) - \frac{1}{\gamma s}\right) \label{Langevin}
\end{equation}
in the limit $\beta_0 \rightarrow 0$, where $M$ and $\chi_0$ are respectively the magnetic saturation and initial susceptibility of the ferrofluid and $\gamma=3\chi_0/M$. Figure \ref{c3 sign} shows that there is a critical magnetic saturation $M_{\mathrm{c}}\approx 10$
below which no localised one-dimensional interfaces bifurcate from the trivial state.

\begin{remark}[\bf Mathematics of the normal-form theory] \label{Math of NFT}
The above choice of coordinates $(A,B,\bar{A},\bar{B})$ is accomplished by the Birkhoff normal-form theory for
$(U_1,\Upsilon,H_\mathrm{f}^\varepsilon)$, which states
that for each $n_0 \geq 2$ there is a near-identity, analytic, symplectic change of coordinates with the property that
$H_\mathrm{f}^\varepsilon$ takes the form \eqref{NF Hamiltonian} in the new coordinates (see Buffoni \& Groves
\cite[pp.\ 196--197]{BuffoniGroves99}). We incorporate this
feature into the construction of $\tilde{r}$ by a replacing the Darboux transformation used in
formula \eqref{Definition of tilde r} by its composition with the normal-form transformation.
\end{remark}

\noindent
{\bf Acknowledgements.} A. S. was supported by the Deutsche Forschungsgemeinschaft under grant GR 3348/1-1. D. J. B. L. acknowledges support from an EPSRC grant (Nucleation of Ferrosolitons and Ferropatterns, EP/H05040X/1). No new data were created during this study. \vspace{-2.5mm}


\end{document}